\documentclass[a4paper, 12pt, reqno]{amsart}

\usepackage{amsmath}
\usepackage{amsthm}
\usepackage{amssymb}
\usepackage{setspace}
\usepackage{enumerate}
\usepackage{comment}
\usepackage{bm}

\newtheorem{thm}{Theorem}[section]
\newtheorem{cor}[thm]{Corollary}
\newtheorem{prop}[thm]{Proposition}
\newtheorem{lem}[thm]{Lemma}

\theoremstyle{definition}

\newtheorem{rem}[thm]{Remark}

\numberwithin{equation}{section}

\newcommand{\pt}{\partial}

\newcommand{\R}{\mathbb{R}}
\newcommand{\C}{\mathbb{C}}

\newcommand{\F}{\mathcal{F}}

\renewcommand{\epsilon}{\varepsilon}

\DeclareMathOperator{\im}{Im}

\DeclareMathOperator*{\weaklim}{w-lim}

\title[1D focusing fourth-order NLS]{Scattering for the focusing, $L^2$-supercritical\\ fourth-order NLS in one dimension}

\author{Koichi Komada}
\address{Department of mechanical and systems engineering,
Chukyo University, Nagoya Aichi 466-8666, Japan}
\email{k-komada@mng.chukyo-u.ac.jp}

\author{Satoshi Masaki}
\address{Department of mathematics, 
Hokkaido University, Sapporo Hokkaido, 060-0810, Japan}
\email{masaki@math.sci.hokudai.ac.jp}



\date{}

\begin{document}

\maketitle

\vskip5mm
\noindent
{\bf Abstract.} In this paper, we consider the fourth-order Schr\"{o}dinger equations with focusing, $L^2$-supercritical nonlinearity in one dimension. We prove the global existence and scattering of  solutions below the ground state threshold under the evenness assumption. This result extends the result by Guo and Dinh in higher dimensions under the radial assumption.

\section{Intrduction}\label{sec:intro}

The classical mathematical model for propagation of intense laser beams in a bulk medium with Kerr nonlinearity is given by the nonlinear Schr\"{o}dinger equation. The role of small fourth-order dispersion has been considered by Karpman \cite{K} and Karpman--Shagalov \cite{KS}. They studied the fourth-order nonlinear Schr\"{o}dinger equations written as
\begin{equation}\label{eq:4NLS}
      i\pt_t u-\Delta^2 u+\epsilon\Delta u=\lambda|u|^{p-1}u,
\end{equation}
where $\epsilon\ge0$, $\lambda\in\R$, $p>1$ and $u=u(t,x)$ is a complex valued function of $(t,x)\in\R\times\R^d$. Fibich--Ilan--Papanicolaou \cite{FIP} investigated the equations \eqref{eq:4NLS} by  numerical and rigorous analysis. Dispersive estimates for the linear fourth-order Schr\"{o}dinger equations have been obtained in 
\cite{BKS}. 

In this article, we study the global dynamics of \eqref{eq:4NLS}.
There is a number of results in this direction. For the defocusing case $\lambda>0$, Pausader \cite{P07} proved global well-posedness and scattering for all radially symmetric initial data in $H^2(\R^d)$ when $d\geq5$ and $1<p\leq 1+8/(d-4)$. Note that the radial assumption in \cite{P07} is not necessary when $p<1+8/(d-4)$ and is removed later in \cite{MXZ11,P09JFA} in the remaining case, $p=1+8/(d-4)$. 
We refer the reader to \cite{PX} for the low dimensional $L^2$-supercritical case: $1\leq d\leq 4$ and $p>1+8/d$. See also \cite{PS} for the high dimensional $L^2$-critical case, i.e., $d\geq5$ and $p=1+8/d$, with $L^2(\R^d)$ data.

Let us turn to the focusing case, $\lambda<0$. 
It is widely known that ground state solutions play an important role in describing global dynamics of nonlinear dispersive equations.
This is also the case for equation \eqref{eq:4NLS}.
Miao--Xu--Zhao \cite{MXZ09} and Pausader \cite{P09DCDS} proved the global well-posedness and the scattering of \eqref{eq:4NLS} in the energy-critical case, i.e., the case where $d\geq5$ and $p=1+8/(d-4)$, for radially symmetric initial data below the ground state threshold. In \cite{MXZ09, P09DCDS}, the proofs are based on concentration compactness argument in the spirit of Kenig--Merle \cite{KM}. This result is extended to the inter-critical case, $d\ge 2$ and $1+8/d<p<1+8/(d-4)_{+}$, by Guo \cite{G}, whose argument also follows the idea in \cite{KM}. Dinh \cite{D} gives another proof for the result of \cite{G} without using concentration compactness argument. The proofs in \cite{D} are based on the idea of Dodson--Murphy \cite{DM}.
We remark that in both \cite{G,D}, the radial symmetry assumption is still used and the one-dimensional case is excluded.

In this paper, we consider the remaining one dimensional case of the focusing inter-critical case with the evenness assumption.
More specifically, our target equation is the
focusing $L^2$-supercritical fourth-order nonlinear Schr\"{o}dinger equation: 
\begin{equation}\label{eq:1DNL4S}
      \left\{
      \begin{split}
            &i\partial_t u-\partial_x^4 u=-|u|^{p-1}u, \ \ \ &&(t,x)\in\R\times\R, \\
            &u(0,x)=u_0(x), \ \ \ &&x\in\R,
      \end{split}
      \right.
\end{equation}
where $u=u(t,x)$ is complex valued, $(t,x)\in\R\times\R$, and $p>9$. 
We concentrate on the case $\varepsilon=0$.
It is known that \eqref{eq:1DNL4S} is locally well-posed in $H^2(\R)$ and solutions to \eqref{eq:1DNL4S} have conserved mass
\begin{align*}
      M(u(t)):=\int_{\R}|u(t,x)|^2 dx
      =M(u_0),
\end{align*}
and energy
\begin{align*}
      E(u(t)):=\frac{1}{2}\int_{\R}|\partial_x^2 u(t,x)|^2 dx-\frac{1}{p+1}\int_{\R}|u(t,x)|^{p+1} dx
      =E(u_0).
\end{align*}
The range $p>9$ corresponds to the mass-supercritical ($L^2$-supercritical) case in the following sense:
The equation \eqref{eq:1DNL4S} is invariant under the scaling transformation
\begin{align*}
      u_{\lambda}(t,x):=\lambda^{\frac{4}{p-1}}u(\lambda^4 t, \lambda x), \ \ \ \lambda>0.
\end{align*}
By direct computation, 
\begin{align*}
      \|u_{\lambda}(0)\|_{\dot{H}^s(\R)}=\lambda^{s-\frac{1}{2}+\frac{4}{p-1}}\|u_0\|_{\dot{H}^s(\R)}
\end{align*}
for all $s\in\R$. This implies that the critical exponent of the homogeneous Sobolev space is given by
\begin{align*}
      s_c:=\frac{1}{2}-\frac{4}{p-1}.
\end{align*}
We say \eqref{eq:1DNL4S} is $L^2$-supercritical if $s_c>0$. The condition reads as $p>9$.

Let us recall the notion of ground states.
\eqref{eq:1DNL4S} has standing wave solutions $u(t,x)=e^{i\omega t}Q_{\omega}(x)$, where $\omega>0$ and $Q_{\omega}$ is a solution of the fourth-order elliptic equation
\begin{align}\label{eq:Q_omega}
      -\omega Q_{\omega}-\pt_x^4 Q_{\omega}+|Q_{\omega}|^{p-1}Q_{\omega}=0,\ \ \ x\in\R.
\end{align}
We say that $Q_{\omega}$ is a ground state of \eqref{eq:Q_omega} if
it is a least-action solution to \eqref{eq:Q_omega}, that is, if it satisfies
the relation
\begin{align*}
      S_{\omega}(Q_{\omega})
      =\inf\{S_{\omega}(f)\ |\ f\in H^2(\R), f\neq0, f\ \mbox{solves \eqref{eq:Q_omega}} \},
\end{align*} 
where $S_{\omega}$ is the action functional defined by
\begin{align*}
      S_{\omega}(f):=E(f)+\frac{\omega}{2}M(f).
\end{align*}
It is known that
there exists a ground state solution to \eqref{eq:Q_omega} for all $p>1$ and $\omega>0$ (see \cite{BCSN,BN}, for instance). 
We remark that if $p$ is an odd integer, then for all $\omega>0$ at least one ground state is even (see \cite{BL}). However, the existence of an even ground state is open for general $p$. 
Let $Q=Q_1$ be one of the ground state solutions to 
\begin{align}\label{eq:Q}
      -Q-\pt_x^4 Q+|Q|^{p-1}Q=0,\ \ \ x\in\R.
\end{align}
Then $Q$ attains the sharp constant $C_{GN}$ of the Gagliardo--Nirenberg inequality
\begin{align}\label{eq:GN}
      \|f\|_{L^{p+1}(\R)}^{p+1}
      \leq C_{GN}\|f\|_{L^2(\R)}^{\frac{3p+5}{4}}\|\pt_x^2 f\|_{L^2(\R)}^{\frac{p-1}{4}}.
\end{align}

\indent
Our main result in this paper is follows:
\begin{thm}[Main result]\label{thm:main}
Let $p>9$, and let $u_0\in H^2(\R)$ be even. If 
\begin{align}\label{eq:below_ground_state}
      M(u_0)^{\frac{2-s_c}{s_c}}E(u_0)
      &<M(Q)^{\frac{2-s_c}{s_c}}E(Q), \ \ \ \mbox{and}\\ \label{eq:below_gs}
      \|u_0\|_{L^2}^{\frac{2-s_c}{s_c}}\|\partial_x^2 u_0\|_{L^2}
      &<\|Q\|_{L^2}^{\frac{2-s_c}{s_c}}\|\partial_x^2 Q\|_{L^2},
\end{align} 
then there exists an unique global solution $u\in C(\R,H^2(\R))$ to \eqref{eq:1DNL4S} with initial data $u_0$. Moreover, the solution $u(t)$ scatters forward and backward in time, i.e., there exist $\phi^{\pm}\in H^2(\R)$ such that 
\begin{align*}
      \lim_{t\rightarrow\pm\infty}\|u(t)-e^{-it\partial_x^4}\phi^{\pm}\|_{H^2}=0.
\end{align*}
\end{thm}

\begin{rem}
\begin{enumerate}
\item
We observe that the scattering conditions \eqref{eq:below_ground_state} and \eqref{eq:below_gs} are sharp for even solutions only when an even ground state exists. Boulenger--Lenzmann \cite{BL} established the existence of radially symmetric ground states for the fourth-order nonlinear Schr\"{o}dinger equations \eqref{eq:4NLS} when $d\geq1$ and $p$ is an odd integer. 
Consequently, our result is sharp in this case.
However, as mentioned earlier, if $p$ is not odd, the existence of such ground states remains unknown due to the inapplicability of classical rearrangement techniques caused by the presence of the fourth-order derivative $\Delta^2$. Their proof relies on the symmetric-decreasing rearrangement in Fourier space. This Fourier rearrangement method is limited to cases where the nonlinearity is a polynomial, making it applicable only when $p$ is an odd number.
\item
We also refer to the symmetry breaking results for ground states of \eqref{eq:4NLS} with $\epsilon<0$. Lenzmann--Weth \cite{LW} proved that all least-action solutions to the fourth-order elliptic equation
\begin{align*}
      -\omega Q_{\omega}-\Delta^2 Q_{\omega}+\epsilon\Delta Q_{\omega}
      +|Q_{\omega}|^{p-1}Q_{\omega}=0,\ \ \ x\in\R^d
\end{align*}
fail to be radially symmetric when $d\geq2$, $\epsilon<0$, $1<p<1+4/(d-1)$ and $\omega$ is slightly larger than $\epsilon^2/4$ (see also \cite{MS}). It is expected that the symmetry breaking for ground states of \eqref{eq:4NLS} does not occur in our case, i.e., when $\epsilon\geq0$. 
\end{enumerate}
\end{rem}

\begin{rem}
It is conjectured that solutions of \eqref{eq:1DNL4S} with $p>9$ blows up in finite time if the initial data $u_0$ satisfies
\begin{align*}
      M(u_0)^{\frac{2-s_c}{s_c}}E(u_0)&\ <M(Q)^{\frac{2-s_c}{s_c}}E(Q),\ \ \ \mbox{and} \\
      \|u_0\|_{L^2}^{\frac{2-s_c}{s_c}}\|\pt_x^2 u_0\|_{L^2}
      &\ >\|Q\|_{L^2}^{\frac{2-s_c}{s_c}}\|\pt_x^2 Q\|_{L^2}.
\end{align*}
It is still an open question. In the case $d\geq2$, Boulenger--Lenzmann \cite{BL} proved finite time blow-up for radially symmetric solutions to \eqref{eq:4NLS} for $1+8/d<p<1+8/(d-4)_{+}$ with $p\leq9$. 
Their proofs are based on localized virial identities for \eqref{eq:4NLS}, which are analogous to virial identities introduced by Ogawa--Tsutsumi \cite{OT} for the second-order nonlinear Schr\"{o}dinger equations. To establish the localized virial identities, the radial Sobolev inequality was employed. Because of this, the radial assumption with $d\geq2$ is required.
\end{rem}

As mentioned above, the radial assumption is used also in the previous studies investigating
the scattering below the ground state threshold for fourth-order nonlinear Schr\"{o}dinger equations.
One motivation for
the use of the radial assumption can be explained by the absence of Galilean invariance, a characteristic property of fourth-order Schr\"{o}dinger equations.
The proof of Theorem \ref{thm:main} is based on the argument by Kenig--Merle \cite{KM},
as in \cite{MXZ09,P09DCDS,G}.
Roughly speaking, we first reformulate the problem as a variational problem and then prove the theorem by showing that the failure of the theorem implies the existence of a phantom optimizer to the problem and that the existence of the phantom yields a contradiction.
The phantom optimizer is often referred to as a \emph{minimal blowup solution} or a \emph{critical element}.

We recall that, in the usual second-order nonlinear Schr\"{o}dinger case,
a similar scattering result is obtained without the radial symmetry for all inter-critical case: $d\ge1$ and $1 + 4/d < p < 1 + 4/(d-2)_+$ (see \cite{AN,FXC,HR,DHR}). 
In \cite{AN,FXC,DHR}, the radial assumption is removed by utilizing the Galilean invariance to assure that the critical element is localized around the origin for all time.
However, in the fourth-order case, the lack of the Galilean invariance makes establishing
the localization of a critical element hard.
The radial assumption immediately stems this step.
This is one advantage of the use of the radial assumption in the previous studies \cite{MXZ09,P09DCDS,G}.

Another benefit of radial symmetry lies in the improvement of the Sobolev embedding.
In both \cite{D} and \cite{G}, the following radial Sobolev inequality was employed: Let $f\in H^1(\R^d)$ be radially symmetric. Then for all $R>0$,
\begin{align}\label{eq:radial_Sobolev}
      \|f\|_{L^{\infty}(|x|>R)}
      \leq CR^{-\frac{d-1}{2}}\|f\|_{L^2(|x|>R)}^{\frac{1}{2}}\|\nabla f\|_{L^2(|x|>R)}^{\frac{1}{2}},
\end{align}
where $C$ is a constant independent of $f$ and $R$. When $d\geq2$, the above estimate \eqref{eq:radial_Sobolev} implies that radially symmetric functions enjoy the stronger spatial decay as $|x|\rightarrow\infty$. 
This improvement makes the analysis simpler in several respects.
For instance, one may deduce that no spatial translation is involved in the linear profile decomposition for sequences of radially symmetric functions. 

The evenness assumption is the one dimensional version 
of the radial assumption in a sense.
However, there is no advantage in view of the uniform decay property such as \eqref{eq:radial_Sobolev}. 
Hence, our main task in this paper is the following two.
One is to modify the concentration compactness argument without using the radial Sobolev inequality \eqref{eq:radial_Sobolev}.
The other is to obtain a long-time localization of a critical element from the evenness assumption.

To achieve the latter, we first capture the improvement of the evenness assumption in the language of the profile decomposition. It then turns out that, although it does not exclude the possibility of the presence of profiles with the space shift, such profiles must appear as a pair (see Remark \ref{rem:improvement}). 
Then, we see that it is strong enough to control the spatial shift of a critical element since it is obtained as an optimizer to a variational problem. An optimizing sequence does not contain more than two profiles. This is the main point of our proof (see Remark \ref{rem:spatialcontrol2}).

The rest of the paper is organized as follows: In Section \ref{sec:preli}, we recall the Strichartz type estimates for the linear fourth-order Schr\"{o}dinger equation. We also give a small data theory and prove a long-time perturbation theory. In Section \ref{sec:variational_analysis}, we study the ground state $Q$ and give two corecivity lemmas. In Section \ref{sec:profile_decomposition}, we prove a linear profile decomposition for sequence of even functions. In Section \ref{sec:concentration_compactness}, we construct a critical solution, which is a target counter example in contradiction argument. Finally, in Section \ref{sec:rigidity}, we prove Theorem \ref{thm:main} by contradiction.

\section{Preliminaries}\label{sec:preli}

We consider the linear fourth-order Schr\"{o}dinger equation
\begin{equation}\label{eq:linear4S}
      \left\{
      \begin{split}
            &i\partial_t u-\partial_x^4 u=0, \ \ \ &&(t,x)\in\R\times\R, \\
            &u(0,x)=u_0(x), \ \ \ &&x\in\R.
      \end{split}
      \right.
\end{equation}
The linear propagator associated with \eqref{eq:linear4S} is given by
\begin{align*}
      e^{-it\partial_x^4}f:=\F^{-1}\left[e^{-it\xi^4}\hat{f}(\xi)\right],
\end{align*}
where $\F f=\hat{f}$ is the Fourier transform of $f$ and $\F^{-1}f$ is the inverse Fourier transform of $f$. Since $|e^{-it\xi^4}|=1$,
\begin{align}\label{eq:L^2}
      \|e^{-it\partial_x^4}f\|_{L^2}=\|f\|_{L^2}.
\end{align}
Ben-Artzi--Koch--Saut \cite{BKS} studied the dispersive estimates for Schr\"{o}dinger equations with fourth-order dispersion. Applying the result in \cite{BKS}, 
\begin{align}\label{eq:L^infty}
      \|e^{-it\partial_x^4}f\|_{L^{\infty}}\lesssim |t|^{-\frac{1}{4}}\|f\|_{L^1}.
\end{align}
Interpolating \eqref{eq:L^2} and \eqref{eq:L^infty}, for $2\leq r\leq\infty$,
\begin{align}\label{eq:dispersive}
      \|e^{-it\partial_x^4}f\|_{L^r}\lesssim |t|^{-\frac{1}{4}(1-\frac{2}{r})}\|f\|_{L^{r^{\prime}}}.
\end{align}

\begin{lem}[Strichartz estimate]\label{lem:Strichartz}
Suppose that pairs $(q_1,r_1)$ and $(q_2,r_2)$ satisfy
\begin{align*}
      2\leq q_j,r_j\leq\infty \ \ \ \mbox{and}\ \ \ \frac{4}{q_j}+\frac{1}{r_j}=\frac{1}{2}, \ \ \ j=1,2.
\end{align*}
Then for any $I\subset\R$ and $t_0\in I$,
\begin{align*}
      &\|e^{-it\partial_x^4}f\|_{L_t^{q_1}L_x^{r_1}(I\times\R)}
      \lesssim \|f\|_{L_x^2(\R)}, \\
      &\left\|\int_{t_0}^t e^{-i(t-s)\partial_x^4}F(s) ds\right\|_{L_t^{q_1}L_x^{r_1}(I\times\R)}
      \lesssim \|F\|_{L_t^{q_2^{\prime}}L_x^{r_2^{\prime}}(I\times\R)}.
\end{align*}
\end{lem}

\noindent
{\bf Proof.} We can prove Lemma \ref{lem:Strichartz} by using \eqref{eq:dispersive} and the theorem in Keel--Tao \cite{KT}.
\qed

\begin{lem}[Kato type estimate]\label{lem:Kato}
Suppose that pairs $(q_1,r_1)$ and $(q_2,r_2)$ satisfy
\begin{align*}
      1\leq q_j<\infty, \ \ \ 2\leq r_j\leq\infty, \ \ \ \frac{2}{q_j}+\frac{1}{r_j}<\frac{1}{2}, \ \ \ j=1,2
\end{align*}
and
\begin{align*}
      \frac{1}{q_1}+\frac{1}{q_2}=\frac{1}{4}\left(1-\frac{1}{r_1}-\frac{1}{r_2}\right).
\end{align*}
Then 
\begin{align*}
      \left\|\int_{t_0}^t e^{-i(t-s)\partial_x^4}F(s) ds\right\|_{L_t^{q_1}L_x^{r_1}(I\times\R)}
      \lesssim \|F\|_{L_t^{q_2^{\prime}}L_x^{r_2^{\prime}}(I\times\R)}.
\end{align*}
\end{lem}

\noindent
{\bf Proof.} We can prove Lemma \ref{lem:Kato} by using \eqref{eq:dispersive} and the theorem in Foschi \cite{F}.
\qed \\

\indent
Throughout this paper, for $p>9$ we fix
\begin{align*}
      &r=p+1, &&q=\frac{8(p+1)}{p-1}, \notag \\
      &q_1=\frac{4(p-1)(p+1)}{3p+5}, &&q_2=\frac{4(p-1)(p+1)}{p^2-5p-4}.
\end{align*}
Then $(q,r)$ satisfies $\frac{4}{q}+\frac{1}{r}=\frac{1}{2}$ and $(q_1,r)$ and $(q_2,r)$ satisfy 
\begin{align*}
      \frac{4}{q_1}+\frac{1}{r}=\frac{4}{p-1}=\frac{1}{2}-s_c,\ \ \ 
      \frac{4}{q_2}+\frac{1}{r}=1-\frac{4}{p-1}=\frac{1}{2}+s_c.
\end{align*}
In particular, we have $\frac{1}{q_1}+\frac{1}{q_2}=\frac{1}{4}(1-\frac{2}{r})$. For any time interval $I\subset\R$, we define function spaces $S(I)$, $X(I)$ with norms
\begin{align*}
      &\|u\|_{S(I)}:= \|u\|_{L_t^{\infty}L_x^2(I\times\R)}+\|u\|_{L_t^qL_x^r(I\times\R)}, \\
      &\|u\|_{X(I)}:= \|u\|_{L_t^{q_1}L_x^r(I\times\R)}.
\end{align*}
Then by the Sobolev embedding, $\|u\|_{X(I)}\lesssim \|D^{s_c}u\|_{S(I)}$. We also define
\begin{align*}
      \|F\|_{N(I)}:=\|F\|_{L_t^{q_2^{\prime}}L_x^{r^{\prime}}(I\times\R)}.
\end{align*}

\begin{lem}[Nonlinear estimate]\label{lem:nonlinear_est}
Let $t_0\in\R$ and let $I$ be a time interval which contains $t_0$. Then
\begin{align*}
      \left\|\int_{t_0}^t e^{-i(t-s)\partial_x^4}(uv)(s) ds\right\|_{X(I)}
      \lesssim \|u\|_{L_t^{q_1}L_x^{r}(I\times\R)}\|v\|_{L_t^{q_1/(p-1)}L_x^{r/(p-1)}(I\times\R)}.
\end{align*}
\end{lem}

\noindent
{\bf Proof.} By using Lemma \ref{lem:Kato},
\begin{align*}
      \left\|\int_{t_0}^t e^{-i(t-s)\partial_x^4}(uv)(s) ds\right\|_{X(I)}
      \lesssim \|uv\|_{N(I)}.
\end{align*}
Since $q_2^{\prime}=\frac{4(p-1)(p+1)}{p(3p+5)}=q_1/p$ and $r^{\prime}=\frac{p+1}{p}=r/p$, by using the H\"{o}lder inequality, we have 
\begin{align*}
      \|uv\|_{N(I)}
      \leq \|u\|_{L_t^{q_1}L_x^r(I\times\R)}\|v\|_{L_t^{q_1/(p-1)}L_x^{r/(p-1)}(I\times\R)}.
\end{align*}
This completes the proof of Lemma \ref{lem:nonlinear_est}.
\qed \\

By using Lemma \ref{lem:nonlinear_est}, we have the following small data theory.

\begin{prop}[Small data]\label{prop:small_data}
Let $p>9$, and let $I$ be a time interval which contains $0$. There exists small $\delta_{sd}>0$ such that if $u_0\in H^2(\R)$ satisfies
\begin{align*}
      \|e^{-it\partial_x^4}u_0\|_{X(I)}\leq\delta_{sd},
\end{align*}
then there exists a unique solution $u\in C(I,H^2(\R))$ to \eqref{eq:1DNL4S} with initial data $u_0$. Moreover,
\begin{align*}
      &\|u\|_{X(I)}\leq 2\|e^{-it\partial_x^4}u_0\|_{X(I)}, \\
      &\|D^{s_c}u\|_{S(I)}\lesssim \|u_0\|_{\dot{H}^{s_c}}.
\end{align*}
\end{prop}

\noindent
{\bf Proof.} We can prove Proposition \ref{prop:small_data} by using Lemma \ref{lem:nonlinear_est} and applying the standard contraction mapping principle. 
\qed \\

In the proof of Theorem \ref{thm:main}, we will use the following scattering criterion.

\begin{prop}[$H^2$-scattering]\label{prop:H^2-scattering}
Let $p>9$. If $u\in C(\R,H^2(\R))$ is a global solution to \eqref{eq:1DNL4S} satisfying
\begin{align}\label{eq:uniform_bound}
      \|u\|_{X(\R)}<\infty, \ \ \ \mbox{and}\ \ \ \|u\|_{L_t^{\infty}H_x^2(\R\times\R)}<\infty,
\end{align}
then $u$ scatters forward and backward in time.
\end{prop}

\noindent
{\bf Proof.} We will prove only that $u$ scatters forward in time. The proof for the backward case is similar. Since $u(t)$ solves the integral equation
\begin{align}\label{eq:integral_eq}
      u(t)=e^{-it\pt_x^4}u(0)+i\int_0^t e^{-i(t-s)\pt_x^4}(|u|^{p-1}u)(s)ds,
\end{align}
we have
\begin{align}\label{eq:integral_eq_scattering}
      u(t)-e^{-it\pt_x^4}\phi^{+}=-i\int_t^{+\infty} e^{-i(t-s)\pt_x^4}(|u|^{p-1}u)(s)ds,
\end{align}
where
\begin{align*}
      \phi^{+}:=u(0)+i\int_0^{+\infty} e^{is\pt_x^4}(|u|^{p-1}u)(s)ds.
\end{align*}
Let $a=\frac{8(p-1)}{7},\ b=2(p-1)$. Then by applying Lemma \ref{lem:Strichartz} to \eqref{eq:integral_eq_scattering}, we obtain
\begin{align}
      \|u(t)-e^{-it\pt_x^4}\phi^{+}\|_{H^2}
      \lesssim&\ \|(1-\Delta)(|u|^{p-1}u)\|_{L_t^{\frac{8}{7}}L_x^1([t,+\infty)\times\R)} \notag \\
      \lesssim&\ \|u\|_{L_t^{\infty}H_x^2}\|u\|_{L_t^a L_x^b([t,+\infty)\times\R)}.
      \label{eq:estimate_H^2_scattering}
\end{align}
Hence, if $\|u\|_{L_t^a L_x^b(\R\times\R)}<\infty$, then \eqref{eq:estimate_H^2_scattering} implies that $\|u(t)-e^{-it\pt_x^4}\phi^{+}\|_{H^2}\rightarrow0$ as $t\rightarrow+\infty$, thus, $u$ scatters forward in time. Since $\frac{4}{a}+\frac{1}{b}=\frac{4}{p-1}=\frac{1}{2}-s_c$, by applying Lemma \ref{lem:Strichartz} and Lemma \ref{lem:Kato} to \eqref{eq:integral_eq} and by using the Sobolev inequality, we have
\begin{align*}
      \|u\|_{L_t^a L_x^b(\R\times\R)}
      \lesssim&\ \|u(0)\|_{\dot{H}^{s_c}}+\||u|^{p-1}u\|_{N(\R)} \notag \\
      \lesssim&\ \|u\|_{L_t^{\infty}H_x^2}+\|u\|_{X(\R)}^p.
\end{align*}
Therefore, from \eqref{eq:uniform_bound} we get $\|u\|_{L_t^a L_x^b(\R\times\R)}<\infty$. This completes the proof of Proposition \ref{prop:H^2-scattering}.
\qed \\

The following long-time perturbation theory plays crucial role in the concentration compactness argument (in Section \ref{sec:concentration_compactness}).

\begin{prop}[Long-time perturbation theory]\label{prop:stability}
Let $p>9$. For each $A\gg1$, there exists $\delta_0=\delta_0(A)\ll1$ and $C=C(A)\gg1$ such that the following holds. Let $I\subset\R$ be a compact time interval, and let $\tilde{u}\in C(I,H^2(\R))$ be an approximate solution of \eqref{eq:1DNL4S} in the sense that
\begin{align*}
      i\pt_t \tilde{u}-\pt_x^4 \tilde{u}=-|\tilde{u}|^{p-1}\tilde{u}+e,
\end{align*}
for some $e\in N(I)$. Let $t_0\in I$ and $u_0\in H^2(\R)$. If 
\begin{align*}
      \|\tilde{u}\|_{X(I)}\leq A,\ \ \ \|e\|_{N(I)}\leq\delta
\end{align*}
and 
\begin{align*}
      \|e^{-i(t-t_0)\pt_x^4}(\tilde{u}(t_0)-u_0)\|_{X(I)}\leq\delta
\end{align*}
for some $\delta\in(0,\delta_0]$, then there exists a solution $u\in C(I,H^2(\R))$ to \eqref{eq:1DNL4S} such that $u(t_0)=u_0$. Moreover, $u$ satisfies
\begin{align*}
      \|u-\tilde{u}\|_{X(I)}\leq C\delta.
\end{align*}
\end{prop}

\noindent
{\bf Proof.} Put $u=\tilde{u}+w$ and let $F(z)=|z|^{p-1}z$. If $w$ solves the equation
\begin{align}\label{eq:w}
      i\pt w-\pt_x^4 w=F(\tilde{u}+w)-F(\tilde{u})-e
\end{align} 
with $w(t_0)=u_0-\tilde{u}(t_0)$, then $u$ solves \eqref{eq:1DNL4S} with $u(t_0)=u_0$. Without loss of generality, we may assume $\inf I=t_0$. Since $\|\tilde{u}\|_{X(I)}\leq A$, we can partition $I$ into $N=N(A)$ intervals $I_j=[t_j,t_{j+1}]$ such that for each $j$,
\begin{align}\label{eq:partition_eta}
      \|\tilde{u}\|_{X(I_j)}\leq\eta
\end{align}
for some $\eta>0$ sufficiently small to be chosen later. The integral equation of \eqref{eq:w} with initial time $t_j$ is given by
\begin{align}\label{eq:integral_w}
      w(t)=e^{-i(t-t_j)\pt_x^4}w(t_j)
      -i\int_{t_j}^t e^{-i(t-s)\pt_x^4}[F(\tilde{u}+w)-F(\tilde{u})-e](s)ds.
\end{align}
We define
\begin{align*}
      \Phi_j(w)(t)
      :=e^{-i(t-t_j)\pt_x^4}w(t_j)
      -i\int_{t_j}^t e^{-i(t-s)\pt_x^4}[F(\tilde{u}+w)-F(\tilde{u})-e](s)ds,
\end{align*}
and let $B_j:=\|e^{-i(t-t_j)\pt_x^4}w(t_j)\|_{X(I)}$. Then by using Lemma \ref{lem:Kato}, we have
\begin{align*}
      \|\Phi_j(w)\|_{X(I_j)}
      \leq&\ \|e^{-i(t-t_j)\pt_x^4}w(t_j)\|_{X(I_j)}
      +c\|F(\tilde{u}+w)-F(\tilde{u})\|_{N(I_j)}+c\|e\|_{N(I_j)} \notag \\
      \leq&\ B_j+c\delta+c\|F(\tilde{u}+w)-F(\tilde{u})\|_{N(I_j)}, 
\end{align*}
where $c$ is a constant depend on $r=p+1$. Let 
\begin{align*}
X_j:=\{ u\in X(I_j)\ |\ \|u\|_{X(I_j)}\leq 2(B_j+c\delta) \}.
\end{align*} 
By Lemma \ref{lem:nonlinear_est} and \eqref{eq:partition_eta},
\begin{align*}
      \|F(\tilde{u}+w)-F(\tilde{u})\|_{N(I_j)}
      \lesssim&\ \left(\|\tilde{u}\|_{X(I_j)}+\|w\|_{X(I_j)}\right)^{p-1}\|w\|_{X(I_j)} \notag \\
      \lesssim&\ \left(\eta+\|w\|_{X(I_j)}\right)^{p-1}\|w\|_{X(I_j)}.
\end{align*}
Hence, if $\eta$ and $B_j+c\delta$ are sufficiently small, for all $w\in X_j$ we have
\begin{align*}
      \|\Phi_j(w)\|_{X(I_j)}\leq 2(B_j+c\delta)
\end{align*}
and so $\Phi_j(w)\in X_j$. Moreover, if $\eta$ and $B_j+c\delta$ sufficiently small, for all $w_1, w_2\in X_j$ we have
\begin{align*}
      \|w_1-w_2\|_{X(I_j)}
      \leq c\|F(\tilde{u}+w_1)-F(\tilde{u}+w_2)\|_{N(I_j)} 
      \leq \frac{1}{2}\|w_1-w_2\|_{X(I_j)}
\end{align*}
and so $\Phi_j: X_j\rightarrow X_j$ is a contraction mapping. Hence, applying the fixed point theorem, there exists an unique $w\in X_j$ such that $\Phi_j(w)=w$. Now take $t=t_{j+1}$ in \eqref{eq:integral_w} and apply $e^{-i(t-t_j)\pt_x^4}$ to both side to obtain
\begin{align*}
      e^{-i(t-t_{j+1})\pt_x^4}w(t_{j+1})
      =&\ e^{-i(t-t_j)\pt_x^4}w(t_j) \notag \\
      &\ -i\int_{t_j}^{t_{j+1}} e^{-i(t-s)\pt_x^4}[F(\tilde{u}+w)-F(\tilde{u})-e](s)ds.
\end{align*}
Since the above Duhamel integral is confined to $I_j=[t_j,t_{j+1}]$, by using Lemma \ref{lem:Kato}, we have
\begin{align*}
      B_{j+1}
      \leq B_j+c\delta+c\|F(\tilde{u}+w)-F(\tilde{u})\|_{N(I_j)}
      \leq 2(B_j+c\delta),
\end{align*}
if $\eta$ and $B_j+c\delta$ are sufficiently small. Iterating this, we obtain
\begin{align*}
      B_j \leq 2^j B_0+c\delta\sum_{k=1}^j 2^k \leq 2^j \delta+2(2^j-1)c\delta
\end{align*}
for all $0\leq j\leq N-1$. Thus, we can take $\delta=\delta(N)>0$ sufficiently small so that by an inductive argument, we get $w\in X_j$ such that $\Phi_j(w)=w$ for all $0\leq j\leq N-1$. Therefore, we obtain a solution $w$ of \eqref{eq:w} such that
\begin{align*}
      \|w\|_{X(I_j)}\leq 2(B_j+c\delta)\leq 2^{j+1}\delta+2(2^{j+1}-1)c\delta
\end{align*}
for all $0\leq j\leq N-1$. Hence,
\begin{align*}
      \|w\|_{X(I)}\leq \delta\sum_{j=0}^{N-1} (2^{j+1}+2(2^{j+1}-1)c).
\end{align*}
This completes the proof of Proposition \ref{prop:stability}.
\qed

\section{Variational analysis}\label{sec:variational_analysis}

Multiplying \eqref{eq:Q} by $Q$, integrating and applying integration by parts, we have
\begin{align}\label{eq:identity}
      \|Q\|_{L^2}^2+\|\pt_x^2Q\|_{L^2}^2-\|Q\|_{L^{p+1}}^{p+1}=0.
\end{align}
Multiplying \eqref{eq:Q} by $x\cdot\nabla Q$, integrating and applying integration by parts, we also have
\begin{align}\label{eq:Pohozhaev}
      -\frac{1}{2}\|Q\|_{L^2}^2-\frac{3}{2}\|\pt_x^2 Q\|_{L^2}^2+\frac{1}{p+1}\|Q\|_{L^{p+1}}^{p+1}=0.
\end{align}
From \eqref{eq:identity} and \eqref{eq:Pohozhaev}, we obtain
\begin{align*}
      \|Q\|_{L^{p+1}}^{p+1}=\frac{4(p+1)}{p-1}\|\pt_x^2 Q\|_{L^2}^2.
\end{align*}
Hence, the sharp constant $C_{GN}$ in \eqref{eq:GN} is given by
\begin{align}\label{eq:C_GN}
      C_{GN}
      =\frac{\|Q\|_{L^{p+1}}^{p+1}}{\|\pt_x^2 Q\|_{L^2}^{\frac{p-1}{4}}\|Q\|_{L^2}^{\frac{3p+5}{4}}}
      =\frac{4(p+1)}{p-1}\left(\|\pt_x^2 Q\|_{L^2}\|Q\|_{L^2}^{\frac{2-s_c}{s_c}}\right)^{-\frac{p-9}{4}}.
\end{align}

\begin{lem}[Corecivity I]\label{lem:corecivity_1}
Let $p>9$. Let $u\in C(I,H^2(\R))$ be the solution to \eqref{eq:1DNL4S} with the maximal time interval $I$ of existence. If 
\begin{align}
      M(u)^{\frac{2-s_c}{s_c}}E(u)
      \leq(1-\delta)M(Q)^{\frac{2-s_c}{s_c}}E(Q), \label{eq:scale_prod._1}
\end{align}
for some $\delta\in(0,1)$ and
\begin{align}
      \|u(t_0)\|_{L^2}^{\frac{2-s_c}{s_c}}\|\pt_x^2 u(t_0)\|_{L^2}
      <\|Q\|_{L^2}^{\frac{2-s_c}{s_c}}\|\pt_x^2 Q\|_{L^2} \label{eq:scale_prod._2}
\end{align}
for some $t_0\in I$, then there exists $\delta^{\prime}=\delta^{\prime}(\delta)\in(0,1)$ such that
\begin{align*}
      \|u(t)\|_{L^2}^{\frac{2-s_c}{s_c}}\|\pt_x^2 u(t)\|_{L^2}
      \leq(1-\delta^{\prime})\|Q\|_{L^2}^{\frac{2-s_c}{s_c}}\|\pt_x^2 Q\|_{L^2}
\end{align*}
for all $t\in I$. In particular, $I=\R$, that is, $u$ exists globally in time.
\end{lem}

\noindent
{\bf Proof.} From the Gagliardo--Nirenberg inequality \eqref{eq:GN} with the sharp constant $C_{GN}$,
\begin{align*}
      M(u)^{\frac{2-s_c}{s_c}}E(u)
      =&\ \frac{1}{2}\|u(t)\|_{L^2}^{2\cdot\frac{2-s_c}{s_c}}\|\pt_x^2 u(t)\|_{L^2}^2
      -\frac{1}{p+1}\|u(t)\|_{L^2}^{2\cdot\frac{2-s_c}{s_c}}\|u(t)\|_{L^{p+1}}^{p+1} \notag \\
      \geq&\ \frac{1}{2}\left(\|u(t)\|_{L^2}^{\frac{2-s_c}{s_c}}\|\pt_x^2 u(t)\|_{L^2}\right)^2
      -\frac{C_{GN}}{p+1}\left(\|u(t)\|_{L^2}^{\frac{2-s_c}{s_c}}\|\pt^2 u(t)\|_{L^2}\right)^{\frac{p-1}{4}}.
\end{align*}
Let $f(x):=\frac{1}{2}x^2-\frac{C_{GN}}{p+1}x^{\frac{p-1}{4}}$ for $x\geq0$. Since
\begin{align*}
      f^{\prime}(x)=x\left(1-\frac{C_{GN}(p-1)}{4(p+1)}x^{\frac{p-9}{4}}\right),
\end{align*}
we see that $f^{\prime}(x)=0$ if $x=0$ or $x=x_1:=\left(\frac{4(p+1)}{C_{GN}(p-1)}\right)^{\frac{4}{p-9}}$, and $f(x)$ is monotone increasing in $0\leq x\leq x_1$. From \eqref{eq:C_GN}, we have $x_1=\|Q\|_{L^2}^{\frac{2-s_c}{s_c}}\|\pt_x^2 Q\|_{L^2}$. From \eqref{eq:scale_prod._1} and the conservation of mass and energy, 
\begin{align*}
      f\left(\|u(t)\|_{L^2}^{\frac{2-s_c}{s_c}}\|\pt_x^2 u(t)\|_{L^2}\right)
      \leq&\ M(u)^{\frac{2-s_c}{s_c}}E(u) \\
      \leq&\ (1-\delta)M(Q)^{\frac{2-s_c}{s_c}}E(Q) \\
      =&\ (1-\delta)f\left(\|Q\|_{L^2}^{\frac{2-s_c}{s_c}}\|\pt_x^2 Q\|_{L^2}\right)
\end{align*}
for all $t\in I$. Hence, from \eqref{eq:scale_prod._2} and the continuity, there exists $\delta^{\prime}=\delta^{\prime}(\delta)\in(0,1)$ such that
\begin{align*}
      \|u(t)\|_{L^2}^{\frac{2-s_c}{s_c}}\|\pt_x^2 u(t)\|_{L^2}
      \leq (1-\delta^{\prime})\|Q\|_{L^2}^{\frac{2-s_c}{s_c}}\|\pt_x^2 Q\|_{L^2}
\end{align*}
for all $t\in I$. In particular, $\|u(t)\|_{H^2}$ is uniformly bounded. Thus, $u$ can be extend to globally exist.
\qed 

\begin{lem}[Corecivity II]\label{lem:corecivity_2}
Let $p>9$. If $\|f\|_{L^2}^{\frac{2-s_c}{s_c}}\|\pt_x^2 f\|_{L^2}\leq(1-\delta^{\prime})\|Q\|_{L^2}^{\frac{2-s_c}{s_c}}\|\pt_x^2 Q\|_{L^2}$ for some $\delta^{\prime}\in(0,1)$, then there exists $\delta^{\prime\prime}\in(0,1)$ such that
\begin{align}\label{eq:lower_K}
      K(f)\geq2\delta^{\prime\prime}\|\pt_x^2 f\|_{L^2}^2,
\end{align}
where
\begin{align*}
      K(f):=\left.\frac{\pt}{\pt\lambda}S_{\omega}(\lambda^{\frac{1}{2}}f(\lambda\cdot))
      \right|_{\lambda=1}
      =2\|\pt_x^2 f\|_{L^2}^2-\frac{p-1}{2(p+1)}\|f\|_{L^{p+1}}^{p+1}.
\end{align*}
Moreover, this implies that
\begin{align}\label{eq:lower_upper_E}
      \frac{p-9}{2(p-1)}\|\pt_x^2 f\|_{L^2}^2
      \leq E(f)
      \leq \frac{1}{2}\|\pt_x^2 f\|_{L^2}^2.
\end{align}
\end{lem}

\noindent
{\bf Proof.} From the sharp Gagliardo--Nirenberg inequality \eqref{eq:GN} and \eqref{eq:C_GN},
\begin{align*}
      M(f)^{\frac{2-s_c}{s_c}}K(f)
      =&\ 2\|f\|_{L^2}^{2\cdot\frac{2-s_c}{s_c}}\|\pt_x^2 f\|_{L^2}^2
      -\frac{p-1}{2(p+1)}\|f\|_{L^2}^{2\cdot\frac{2-s_c}{s_c}}\|f\|_{L^{p+1}}^{p+1} \notag \\
      \geq&\ 2\left(\|f\|_{L^2}^{\frac{2-s_c}{s_c}}\|\pt_x^2 f\|_{L^2}\right)^2
      -\frac{C_{GN}(p-1)}{2(p+1)}
      \left(\|f\|_{L^2}^{\frac{2-s_c}{s_c}}\|\pt_x^2 f\|_{L^2}\right)^{\frac{p-1}{4}}.
\end{align*}
Let $g(x):=2x^2-\frac{C_{GN}(p-1)}{2(p+1)}x^{\frac{p-1}{4}}$ for $x\geq0$. Then $g(x)=2x^2h(x)$, where $h(x)=1-\frac{C_{GN}(p-1)}{4(p+1)}x^{\frac{p-9}{4}}$. Using \eqref{eq:C_GN}, for $0\leq x\leq(1-\delta^{\prime})\|Q\|_{L^2}^{\frac{2-s_c}{s_c}}\|\pt_x^2 Q\|_{L^2}$,  we have
\begin{align*}
      h(x)=1-\left(\frac{x}{\|Q\|_{L^2}^{\frac{2-s_c}{s_c}}\|\pt_x^2 Q\|_{L^2}}\right)^{\frac{p-9}{4}}
      \geq1-(1-\delta^{\prime})^{\frac{p-9}{4}}.
\end{align*}
Put $\delta^{\prime\prime}=1-(1-\delta^{\prime})^{\frac{p-9}{4}}$. Then $\delta^{\prime\prime}\in(0,1)$ and $g(x)\geq2\delta^{\prime\prime}x^2$. Therefore, if $\|f\|_{L^2}^{\frac{2-s_c}{s_c}}\|\pt_x^2 f\|_{L^2}\leq(1-\delta^{\prime})\|Q\|_{L^2}^{\frac{2-s_c}{s_c}}\|\pt_x^2 Q\|_{L^2}$, then
\begin{align*}
      M(f)^{\frac{2-s_c}{s_c}}K(f)
      \geq2\delta^{\prime\prime}\left(\|f\|_{L^2}^{\frac{2-s_c}{s_c}}\|\pt_x^2 f\|_{L^2}\right)^2.
\end{align*}
This implies \eqref{eq:lower_K}. Then since
\begin{align*}
      E(f)=\frac{p-9}{2(p-1)}\|\pt_x^2 f\|_{L^2}^2+\frac{2}{p-1}K(f),
\end{align*}
we obtain \eqref{eq:lower_upper_E}.
\qed

\section{Profile decomposition}\label{sec:profile_decomposition}

In this section, we derive a linear profile decomposition for the biharmonic Schr\"{o}dinger propagator in one spacial dimension. 

\begin{prop}[Profile decomposition]\label{prop:profile_decomp.}
Let $p>9$. Let $\{f_n\}$ be a sequence of functions uniformly bounded in $H^2(\R)$. Passing to a subsequence if necessary, there exists $J^{*}\in\{0,1,...\}\cup\{\infty\}$ and for each $1\leq j\leq J^{*}$, there exists a function $\phi^j\in H^2(\R)\setminus\{0\}$ and a sequence $(t_n^j,x_n^j)\in\R\times\R$ such that for each finite $1\leq J\leq J^{*}$, we have the decomposition
\begin{align}\label{eq:decomp.}
      f_n(x)=\sum_{j=1}^J e^{it_n^j\pt_x^4}\phi^j(x-x_n^j)+R_n^J(x)
\end{align}
with the following properties: 
\begin{enumerate}\renewcommand{\labelenumi}{(\roman{enumi})}
\item 
\begin{align}\label{eq:smallness}
      \lim_{J\rightarrow J^{*}}\limsup_{n\rightarrow\infty}\|e^{-it\pt_x^4}R_n^J\|_{X(\R)}=0.
\end{align}
\item For each finite $1\leq J\leq J^{*}$,
\begin{align}\label{eq:weakly_to_0}
      e^{-it_n^J\pt_x^4}R_n^J(x+x_n^J)\rightharpoonup0\ \ \ \mbox{weakly in}\ H^2(\R).
\end{align}
\item For each finite $1\leq J\leq J^{*}$,
\begin{align}\label{eq:norm_decomp.}
      \lim_{n\rightarrow\infty}\left\{
      \|f_n\|_{\dot{H}^s}^2-\sum_{j=1}^J \|\phi^j\|_{\dot{H}^s}^2-\|R_n^J\|_{\dot{H}^s}^2
      \right\}=0,
      \ \ \ \mbox{for all}\ 0\leq s\leq2.
\end{align}
\item For each finite $j \in [1,J^*]$, either $x_n^j \equiv 0$ or $\lim_{n\to\infty} |x_n^j|=\infty$ holds. Similarly, for each finite $j \in [1,J^*]$, either $t_n^j \equiv 0$, $\lim_{n\to\infty} t_n^j=\infty$, or $\lim_{n\to\infty} t_n^j=-\infty$ holds.
If $1\leq j\neq k\leq J^{*}$ then
\begin{align}\label{eq:orthogonality}
      |x_n^j-x_n^k|+|t_n^j-t_n^k|\rightarrow\infty \ \ \ \mbox{as}\ \ \ n\rightarrow\infty.
\end{align}
\item Furthermore, if $\{f_n\}$ is a sequence of even functions then the following is valid:
Pick finite $ j\in [1, J^{*}]$. 
If $\lim_{n\to\infty}|x_n^j|=\infty$, then there exists a unique number $ j^{\prime}\in [1, J^{*}]\setminus\{ j\}$ such that $|x_n^j-(-x_n^{j^{\prime}})|+|t_n^j-t_n^{j'}|$ is uniformly bounded.
If $x_n^j \equiv0$ then $\phi^j$ is even.
\end{enumerate}
\end{prop}

\begin{rem}
If $\{f_n\}\subset H^2(\R^d)$ is uniformly bounded and $d\geq2$, then we can obtain a  decomposition similar to \eqref{eq:decomp.} (with obvious modifications such as replacing $\pt_x^4$ by $\Delta^2$). If, in addition, $\{f_n\}\subset H^2(\R^d)$ is a sequence of radially symmetric functions then we see that, for all $j \in [1,J^*]$, there is no space translation, i.e., $x_n^j\equiv0$, and
the profile function $\phi^j$ is radial. This additional property can be proved by using the radial Sobolev inequality \eqref{eq:radial_Sobolev}. However, in contrast, the property $x_n^j\equiv0$ is not obtained in one dimensional case even if $\{f_n\}\subset H^2(\R)$ is a sequence of even functions.
\end{rem}

\begin{rem}\label{rem:improvement}
The last property (v) in Proposition \ref{prop:profile_decomp.} implies that if $\{f_n\}$ is a sequence of even functions, then for each profile $\phi^j$ with unbounded $x_n^j$ there exists another profile $\phi^{j^{\prime}}$ such that $x_n^{j^{\prime}}=-x_n^j+y_n$ and $t_n^{j'} = t_n^j + \tau_n$ hold,
where $y_n$ and $\tau_n$ are uniformly bounded. 
By taking a subsequence and modifying the definition of the profile and the remainder, one may assume that $y_n\equiv \tau_n \equiv 0$.
This improvement plays an important role in our argument.
\end{rem}

\begin{lem}\label{lem:t_n_x_n}
Let $f\in H^2(\R)$, and let $\{(t_n, x_n)\}\subset\R\times\R$ be such that $|t_n|\rightarrow\infty$ or $|x_n|\rightarrow\infty$. Then,
\begin{align*}
      e^{-it_n\pt_x^4}f(x+x_n)\rightharpoonup 0 \ \ \ \mbox{weakly in}\ H^2(\R).
\end{align*}

\end{lem}

\noindent
{\bf Proof.} Lemma \ref{lem:t_n_x_n} can be proved by using a density argument and \eqref{eq:dispersive}. We omit the detail. 
\qed 

\begin{lem}\label{lem:inverse}
Let $p>9$. Let $\{f_n\}$ be a sequence of functions in $H^2(\R)$. Suppose that
\begin{align}\label{eq:A_and_epsilon}
      A:=\lim_{n\rightarrow\infty}\|f_n\|_{H^2}<\infty
      \ \ \ \mbox{and}\ \ \
      \epsilon:=\lim_{n\rightarrow\infty}\|e^{-it\pt_x^4}f_n\|_{L_t^{\infty}L_x^{\frac{p-1}{4}}}>0.
\end{align}
Then, passing to a subsequence, there exist a function $\phi\in H^2(\R)$ and a sequence $\{(t_n,x_n)\}\subset\R\times\R$ such that
\begin{align}
      e^{-it_n\pt_x^4}f_n(x+x_n)\rightharpoonup&\ \phi \ \ \ \mbox{weakly in}\ H^2(\R),
      \label{eq:weaklimit} \\
      \|\phi\|_{\dot{H}^{s_c}}\gtrsim&\ A\left(\frac{\epsilon}{A}\right)^{\frac{3p-11}{p-9}}. \notag
\end{align}

\end{lem}

\noindent
{\bf Proof.} \eqref{eq:A_and_epsilon} implies that for sufficiently large $n$, 
\begin{align*}
      \|f_n\|_{H^2}\leq2A\ \ \ \mbox{and}\ \ \ 
      \|e^{it\pt_x^4}f_n\|_{L_t^{\infty}L_x^{\frac{p-1}{4}}}\geq\frac{\epsilon}{2}.
\end{align*}
Let $r\geq1$ to be chosen later, and let $\chi(x)$ be a smooth radial function such that $\hat{\chi}(\xi)=1$ for $\frac{1}{r}\leq|\xi|\leq r$ and $\hat{\chi}(\xi)$ is supported in $\frac{1}{2r}\leq|\xi|\leq2r$. By the Sobolev embedding,  
\begin{align*}
      \|e^{-it\pt_x^4}f_n-\chi\ast e^{-it\pt_x^4}f_n\|_{L_t^{\infty}L_x^{\frac{p-1}{4}}}^2
      \lesssim&\int_{\R} |\xi|^{2s_c}(1-\hat{\chi}(\xi))^2|\widehat{f_n}(\xi)|^2 d\xi \notag \\
      \leq&\int_{|\xi|\leq\frac{1}{r}}|\xi|^{2s_c}|\widehat{f_n}(\xi)|^2 d\xi
      +\int_{|\xi|\geq\frac{1}{r}}|\xi|^{2s_c}|\widehat{f_n}(\xi)|^2 d\xi \notag \\
      \leq&\ r^{-2s_c}\|f_n\|_{L^2}^2+r^{-2(2-s_c)}\|f_n\|_{\dot{H}^2}^2 .
\end{align*}
Since $r>1$ and $s_c\in(0,1)$, for sufficiently large $n$, we have
\begin{align*}
      \|e^{-it\pt_x^4}f_n-\chi\ast e^{-it\pt_x^4}f_n\|_{L_t^{\infty}L_x^{\frac{p-1}{4}}}
      \leq\ cr^{-s_c}A.
\end{align*}
Thus, taking $r=\max\{1, \left(\frac{4cA}{\epsilon}\right)^{\frac{1}{s_c}}\}$ so that $cr^{-s_c}A\leq\frac{\epsilon}{4}$, 
\begin{align*}
      \|\chi\ast e^{-it\pt_x^4}f_n\|_{L_t^{\infty}L_x^{\frac{p-1}{4}}}.
      \geq \frac{\epsilon}{4}.
\end{align*}
On the other hand, by H\"{o}lder inequality, we obtain
\begin{align*}
      \|\chi\ast e^{-it\pt_x^4}f_n\|_{L_t^{\infty}L_x^{\frac{p-1}{4}}}
      \leq&\ \|\chi\ast e^{-it\pt_x^4}f_n\|_{L_t^{\infty}L_x^2}^{\frac{8}{p-1}}
      \|\chi\ast e^{-it\pt_x^4}f_n\|_{L_{t,x}^{\infty}}^{\frac{p-9}{p-1}} \notag \\
      \leq&\ A^{\frac{8}{p-1}}\|\chi\ast e^{-it\pt_x^4}f_n\|_{L_{t,x}^{\infty}}^{\frac{p-9}{p-1}}.
\end{align*}
Therefore, we have
\begin{align*}
      \|\chi\ast e^{-it\pt_x^4}f_n\|_{L_{t,x}^{\infty}}
      \gtrsim A\left(\frac{\epsilon}{A}\right)^{\frac{p-1}{p-9}}.
\end{align*}
This implies that there exists $(t_n, x_n)\in\R\times\R$ such that
\begin{align}\label{eq:lower_phi_n}
      |\chi\ast e^{-it_n\pt_x^4}f_n(x_n^{J+1})|
      \gtrsim A\left(\frac{\epsilon}{A}\right)^{\frac{p-1}{p-9}}.
\end{align}
Let $g_n:=e^{-it_n\pt_x^4}f_n(x+x_n)$. Then $\|g_n\|_{H^2}=\|f_n\|_{H^2}\leq 2A$. Hence, there exists $\phi\in H^2(\R)$ and $g_n\rightharpoonup\phi$ weakly in $H^2(\R)$. By the H\"{o}lder inequality and the Sobolev inequality and using \eqref{eq:lower_phi_n}, we have
\begin{align}\label{eq:lower_chiphi}
      \|\chi\|_{L^{\frac{p-1}{p-5}}}\|\phi\|_{\dot{H}^{s_c}}
      \gtrsim|\langle\chi, \phi\rangle_{L^2}|
      =&\lim_{n\rightarrow\infty}|\langle\chi, g_n\rangle_{L^2}| \notag \\
      =&\lim_{n\rightarrow\infty}\left|\int_{\R} \chi(x)e^{-it_n\pt_x^4}f_n(x+x_n) dx\right|
      \notag \\
      =&\lim_{n\rightarrow\infty}|\chi\ast e^{-it_n\pt_x^4}f_n(x_n)| \notag \\
      \gtrsim&\ A\left(\frac{\epsilon}{A}\right)^{\frac{p-1}{p-9}}.
\end{align}
Since $\|\chi\|_{L^{\frac{p-1}{p-5}}}\lesssim r^{\frac{p-5}{p-1}}\lesssim \left(\frac{A}{\epsilon}\right)^{\frac{2(p-5)}{p-9}}$, from \eqref{eq:lower_chiphi} we have
\begin{align*}
      \|\phi\|_{\dot{H}^{s_c}}\gtrsim A\left(\frac{\epsilon}{A}\right)^{\frac{3p-11}{p-9}}.
\end{align*}
This completes the proof of Lemma \ref{lem:inverse}.
\qed \\

\noindent
{\bf Proof of Proposition \ref{prop:profile_decomp.}.} For a uniformly bounded sequence $\{f_n\}\subset H^2(\R)$, we define
\begin{equation}\label{eq:set_of_bubbles}
      \mathcal{V}(\{f_n\})
      :=\left\{
            \phi\in H^2(\R)\ \left|\ 
            \begin{split}
                  &\exists\{(t_n, x_n)\}\subset\R\times\R\ \mbox{s.t.}\ 
                  \mbox{passing to a subsequence,} \\
                  &e^{-it_n\pt_x^4}f_n(x+x_n)\rightharpoonup \phi\ \mbox{weakly in}\ H^2(\R).
            \end{split}
            \right.
      \right\}
\end{equation}
and
\begin{align}\label{eq:eta}
      \eta(\{f_n\})
      :=\sup_{\phi\in\mathcal{V}(\{f_n\})}\|\phi\|_{H^2}.
\end{align}
Then, we note that $\eta(\{f_n\})\leq\limsup_{n\rightarrow\infty}\|f_n\|_{H^2}$. To prove Proposition \ref{prop:profile_decomp.}, we first show the decomposition \eqref{eq:decomp.} with the property (i) replaced by
\begin{align}\label{eq:small_eta}
      \lim_{J\rightarrow J^{*}}\eta(\{R_n^J\})=0,
\end{align}
and with properties (ii), (iii) and (iv). We will extract $\phi^j,t_n^j,x_n^j$ by an inductive process. \\
\indent
Set $R_n^0:=f_n$. If $\eta(\{f_n\})=0$, then the proof is complete with $J^{*}=0$. If not, by the definition of $\eta(\{f_n\})$, there exists $\phi^1\in\mathcal{V}(\{f_n\})$ such that 
\begin{align*}
      \|\phi^1\|_{H^2}\geq\frac{1}{2}\eta(\{f_n\})>0.
\end{align*}
In particular, $\phi^1\neq0$. Moreover, by the definition of $\mathcal{V}(\{f_n\})$, there exists $\{(t_n^1,x_n^1)\}\subset\R\times\R$ such that
\begin{align}\label{eq:pdpf1}
      e^{-it_n^1\pt_x^4}f_n(x+x_n^1)\rightharpoonup\phi^1\ \ \ \mbox{weakly in}\ H^2(\R).
\end{align}
Passing to a subsequence, one may suppose that either $\lim_{n\to\infty}|x_n^1|=\infty$ or $x_n^1\to x_1\in \R$ as $n\to\infty$.
We claim that one can choose $x_n^1\equiv 0$ in the latter case. This is because
\[
	\langle e^{-it_n^1\pt_x^4}f_n, \psi \rangle_{H^2}
	= \langle e^{-it_n^1\pt_x^4}f_n(\cdot+x_n^1), \psi(\cdot+x_n^1) \rangle_{H^2} \to
	\langle \phi^1, \psi(\cdot+x^1) \rangle_{H^2}
	= \langle \phi^1(\cdot-x^1), \psi \rangle_{H^2}
\]
as $n\to\infty$. This implies that \eqref{eq:pdpf1} is valid with the choice $x_n^1 \equiv 0$ with the modification $\phi^1 \mapsto \phi^1(\cdot - x^1)$.
The claim is proved. Similarly, one verifies that either $\lim_{n\to\infty}|t_n^1|=\infty$ or $t_n^1 \equiv 0$ holds, by passing to a subsequence if necessary.

We set $R_n^1:=f_n-e^{it_n^1\pt_x^4}\phi^1(x-x_n^1)$. Then we have
\begin{align*}
      e^{-it_n^1\pt_x^4}R_n^1(x+x_n^1)\rightharpoonup0
\end{align*}
and 
\begin{align*}
      \|f_n\|_{\dot{H}^s}^2-\|R_n^1\|_{\dot{H}^s}^2
      =&\ 2\langle f_n, e^{it_n^1\pt_x^4}\phi^1(x-x_n^1)\rangle_{\dot{H}^s}
      -\|e^{it_n^1\pt_x^4}\phi^1(x-x_n^1)\|_{\dot{H}^s}^2 \notag \\
      =&\ 2\langle e^{-it_n^1\pt_x^4}f_n(x+x_n^1), 
      \phi^1\rangle_{\dot{H}^s}-\|\phi^1\|_{\dot{H}^s}^2 \notag \\
      \rightarrow&\ \|\phi^1\|_{\dot{H}^s}^2,
\end{align*}
for any $0\leq s\leq2$. Hence, we have \eqref{eq:weakly_to_0} and \eqref{eq:norm_decomp.} for $J=1$. \\
\indent
Let $J_0\geq1$ be fixed. We suppose that for all  $1\leq J\leq J_0$ we have the decomposition \eqref{eq:decomp.} with \eqref{eq:weakly_to_0} and \eqref{eq:norm_decomp.}, and also suppose that  \eqref{eq:orthogonality} holds for any $1\leq j\neq k\leq J_0$. If $\eta(\{R_n^{J_0}\})=0$, then we stop the induction and set $J^{*}=J_0$. If $\eta(\{R_n^{J_0}\})>0$, then by the definition of $\eta(\{R_n^{J_0}\})$, we find $\phi^{J_0+1}\in\mathcal{V}(\{R_n^{J_0}\})$ such that
\begin{align}\label{eq:nonzero_phi}
      \|\phi^{J_0+1}\|_{H^2}\geq\frac{1}{2}\eta(\{R_n^{J_0}\})>0.
\end{align}
Moreover, there exist $\{(t_n^{J_0+1}, x_n^{J_0+1})\}\subset\R\times\R$ such that
\begin{align*}
      e^{-it_n^{J_0+1}\pt_x^4}R_n^{J_0+1}(x+x_n^{J_0+1})\rightharpoonup\phi^{J_0+1}
      \ \ \ \mbox{weakly in}\ H^2(\R).
\end{align*}
Arguing as above, passing to a subsequence if necessary, we see that $\{t_n^{J_0+1}\}$ and
$\{x_n^{J_0+1}\}$ satisfy the alternative in the property (iv).
We set $R_n^{J_0+1}:=R_n^{J_0}-e^{it_n^{J_0+1}\pt_x^4}\phi^{J_0+1}(x-x_n^{J_0+1})$. Then we have
\begin{align*}
      e^{-it_n^{J_0+1}\pt_x^4}R_n^{J_0+1}(x+x_n^{J_0+1})\rightharpoonup0
\end{align*}
and 
\begin{align}\label{eq:norm_induct.}
      \|R_n^{J_0}\|_{\dot{H}^s}^2-\|R_n^{J_0+1}\|_{\dot{H}^s}^2
      \rightarrow\|\phi^{J_0+1}\|_{\dot{H}^s}^2,
\end{align}
for any $0\leq s\leq2$. Hence, we have the decomposition \eqref{eq:decomp.} with properties \eqref{eq:weakly_to_0} and \eqref{eq:norm_decomp.} for $J=J_0+1$. \\
\indent
To prove \eqref{eq:orthogonality} for all $1\leq j\neq k\leq J_0+1$, it suffices to claim that for all $1\leq j\leq J_0$, 
\begin{align}\label{eq:claim_orthogonality}
      |x_n^j-x_n^{J_0+1}|+|t_n^j-t_n^{J_0+1}|\rightarrow\infty.
\end{align}
Suppose that there exists $1\leq j_0\leq J_0$ such that $|x_n^{j_0}-x_n^{J_0+1}|$ and $|t_n^{j_0}-t_n^{J_0+1}|$ are finite. Then passing to a subsequence, we may assume that 
\begin{align}\label{eq:ifnot_orthogonality}
      x_n^{j_0}-x_n^{J_0+1}\rightarrow x_0
      \ \ \ \mbox{and}\ \ \ 
      t_n^{j_0}-t_n^{J_0+1}\rightarrow t_0.
\end{align}
Since $e^{-it_n^{j_0}\pt_x^4}R_n^{j_0}(x+x_n^{j_0})\rightharpoonup0$, from \eqref{eq:ifnot_orthogonality} we have
\begin{align*}
      &\weaklim_{n\rightarrow\infty}e^{-it_n^{J_0+1}\pt_x^4}R_n^{j_0}(x+x_n^{J_0+1}) \notag \\
      =&\ \weaklim_{n\rightarrow\infty}
      e^{i(t_n^{j_0}-t_n^{J_0+1})\pt_x^4}\left[
            e^{-it_n^{j_0}\pt_x^4}R_n^{j_0}(\cdot+x_n^{j_0})
      \right](x-(x_n^{j_0}-x_n^{J_0+1})) \notag \\
      =&\ 0.
\end{align*}
Hence, if $j_0=J_0$, then $\phi^{J_0+1}=0$. If $J_0\geq2$ and $1\leq j_0\leq J_0-1$, using the inductive relation $R_n^{J_0}=R_n^{j_0}-\sum_{j=j_0+1}^{J_0} e^{it_n^j\pt_x^4}\phi^j(x-x_n^j)$, we have
\begin{align*}
      \phi^{J_0+1}
      =&\weaklim_{n\rightarrow\infty}e^{-it_n^{J_0+1}\pt_x^4}R_n^{J_0}(x+x_n^{J_0+1}) \notag \\
      =&\weaklim_{n\rightarrow\infty}e^{-it_n^{J_0+1}\pt_x^4}R_n^{j_0}(x+x_n^{J_0+1}) \notag \\
      &-\sum_{j=j_0+1}^{J_0} \weaklim_{n\rightarrow\infty}
      e^{i(t_n^j-t_n^{J_0+1})\pt_x^4}\phi^j(x-(x_n^j-x_n^{J_0+1})) \notag \\
      =&-\sum_{j=j_0+1}^{J_0} \weaklim_{n\rightarrow\infty}
      e^{i(t_n^j-t_n^{J_0+1})\pt_x^4}\phi^j(x-(x_n^j-x_n^{J_0+1})).
\end{align*}
Since $|x_n^j-x_n^{j_0}|+|t_n^j-t_n^{j_0}|\rightarrow\infty$ for all $j_0+1\leq j\leq J_0$, from \eqref{eq:ifnot_orthogonality} and Lemma \ref{lem:t_n_x_n},
\begin{align*}
      &\weaklim_{n\rightarrow\infty}e^{i(t_n^j-t_n^{J_0+1})\pt_x^4}\phi^j(x-(x_n^j-x_n^{J_0+1})) 
      \notag \\
      =&\ \weaklim_{n\rightarrow\infty}e^{i(t_n^{j_0}-t_n^{J_0+1})\pt_x^4}\left[
            e^{i(t_n^j-t_n^{j_0})\pt_x^4}\phi^j(\cdot-(x_n^j-x_n^{j_0}))
      \right](x-(x_n^{j_0}-x_n^{J_0+1})) \notag \\
      =&\ 0
\end{align*}
for all $j_0+1\leq j\leq J_0$. Hence, $\phi^{J_0+1}=0$. This contradicts \eqref{eq:nonzero_phi} and so we have proved that \eqref{eq:orthogonality} holds for all $1\leq j\neq k\leq J_0+1$. \\
\indent
If $\eta(\{R_n^{J_0+1}\})=0$, we stop the induction and set $J^{*}=J_0+1$. In this case, \eqref{eq:small_eta} holds automatically. If $\eta(\{R_n^{J_0+1}\})>0$, we continue the induction. If the inductive process does not terminate in finitely many steps, we set $J^{*}=\infty$. In this case, since $\|\phi^{J+1}\|_{H^2}\geq\frac{1}{2}\eta(\{R_n^J\})$ and $\|R_n^J\|_{H^2}^2-\|R_n^{J+1}\|_{H^2}^2\rightarrow\|\phi^{J+1}\|_{H^2}^2$ for each $J\geq0$ in the process, we obtain
\begin{align*}
      \|R_n^{J+1}\|_{H^2}^2
      =&\ \|R_n^J\|_{H^2}^2-\|\phi_n^{J+1}\|_{H^2}^2+o(1) \notag \\
      \leq&\ \frac{3}{4}\eta(\{R_n^J\})^2+o(1),
\end{align*}
where $o(1)\rightarrow0$ as $n\rightarrow\infty$ for each $J\geq0$. Hence, for each $J\geq0$, we have
\begin{align}\label{eq:A_J+1<A_J}
      \eta(\{R_n^{J+1}\})
      \leq \limsup_{n\rightarrow\infty}\|R_n^{J+1}\|_{H^2}
      \leq \frac{\sqrt{3}}{2}\eta(\{R_n^J\}).
\end{align}
This combined with $\eta(\{R_n^0\})=\eta(\{f_n\})\leq\limsup_{n\rightarrow\infty}\|f_n\|_{H^2}<\infty$ implies that $\eta(\{R_n^J\})\rightarrow0$ as $J\rightarrow\infty$ and so \eqref{eq:small_eta} follows. \\
\indent
Now we will prove the property (i). 
Note by the interpolation, the Strichartz estimate (Lemma \ref{lem:Strichartz}) and the Sobolev inequality, it suffices to show
\begin{align}\label{eq:smallness_2}
      \lim_{J\rightarrow J^{*}}\limsup_{n\rightarrow\infty}
      \|e^{-it\pt_x^4}R_n^J\|_{L_t^{\infty}L_x^{\frac{p-1}{4}}(\R\times\R)}=0.
\end{align}
For each finite $0\leq J\leq J^{*}$, let 
\begin{align*}
      A_J:=\limsup_{n\rightarrow\infty}\|R_n^J\|_{H^2}
      \ \ \ \mbox{and}\ \ \ 
      \epsilon_J:=\limsup_{n\rightarrow\infty}
      \|e^{-it\pt_x^4}R_n^J\|_{L_t^{\infty}L_x^{\frac{p-1}{4}}(\R\times\R)}.
\end{align*}
Applying Lemma \ref{lem:inverse}, there exist $\tilde{\phi}^{J+1}\in H^2(\R)$ and $\{(\tilde{t}_n^{J+1}, \tilde{x}_n^{J+1})\}\subset\R\times\R$ such that
\begin{align*}
      e^{-i\tilde{t}_n^{J+1}\pt_x^4}R_n^J(x+\tilde{x}_n^{J+1})\rightharpoonup\tilde{\phi}^{J+1}
      \ \ \ \mbox{weakly in}\ H^2(\R)
\end{align*}
and
\begin{align}\label{eq:lower_phi^J+1}
      \|\tilde{\phi}^{J+1}\|_{\dot{H}^{s_c}}
      \gtrsim A_J\left(\frac{\epsilon_J}{A_J}\right)^{\frac{3p-11}{p-9}}
      =A_J^{-\frac{2(p-1)}{p-9}}\epsilon_J^{\frac{3p-11}{p-9}}.
\end{align}
Since $\tilde{\phi}^{J+1}\in\mathcal{V}(\{R_n^J\})$, by using \eqref{eq:small_eta},
\begin{align}\label{eq:upper_phi^J+1}
      \|\tilde{\phi}^{J+1}\|_{\dot{H}^{s_c}}
      \lesssim \|\tilde{\phi}^{J+1}\|_{H^2}
      \leq \eta(\{R_n^J\})
      \rightarrow0\ \ \ \mbox{as}\ J\rightarrow J^{*}.
\end{align}
From \eqref{eq:A_J+1<A_J}, we see $A_{J+1}\leq A_J$. In particular, $A_J\leq A_0=\limsup_{n\rightarrow\infty}\|f_n\|_{H^2}<\infty$ for all finite $0\leq J\leq J^{*}$. Hence, \eqref{eq:lower_phi^J+1} and \eqref{eq:upper_phi^J+1} imply $\epsilon_J\rightarrow0$ as $J\rightarrow J^{*}$ and so \eqref{eq:smallness_2} follows.

Finally, we will show the improvement under the evenness assumption, 
which corresponds to (v). 
 Suppose that $\{f_n\}$ is a sequence of even functions and that there exists a finite $1\leq j_0\leq J^{*}$ such that $|x_n^{j_0}|\to \infty$ ($n\to\infty$) but
$|x_n^{j_0}-(-x_n^j)|+|t_n^{j_0}-t_n^j|$ is unbounded for all $1\leq j\leq J^{*}$ with $j \neq j_0$. Then passing to a subsequence if necessary, we may assume that $|x_n^{j_0}+x_n^j|\rightarrow\infty$ or $|t_n^{j_0}-t_n^j|\rightarrow\infty$ holds as $n\to\infty$. We first claim
\begin{align}\label{eq:weakly_to_phi^j_0}
      e^{-it_n^{j_0}\pt_x^4}f_n(x+x_n^{j_0})\rightharpoonup\phi^{j_0}(x)
      \ \ \ \mbox{weakly in}\ H^2(\R).
\end{align}
Indeed, by using the properties (ii), (iv) and applying Lemma \ref{lem:t_n_x_n},
\begin{align*}
      e^{-it_n^{j_0}\pt_x^4}f_n(x+x_n^{j_0})
      =&\ \sum_{j=1}^{j_0} e^{i(t_n^j-t_n^{j_0})\pt_x^4}\phi^j(x-(x_n^j-x_n^{j_0}))
      +e^{-it_n^{j_0}\pt_x^4}R_n^{j_0}(x+x_n^{j_0}) \notag \\
      \rightharpoonup&\ \phi^{j_0}(x).
\end{align*}
Since $f_n$ is even, from \eqref{eq:weakly_to_phi^j_0} we have
\begin{align}\label{eq:pdevenpf}
      e^{-it_n^{j_0}\pt_x^4}f_n(x+(-x_n^{j_0}))
      =&\ e^{-it_n^{j_0}\pt_x^4}f_n(-x+x_n^{j_0}) \notag \\
      \rightharpoonup&\ \phi^{j_0}(-x).
\end{align}
On the other hand, for any finite $1\leq J\leq J^{*}$, 
\begin{align*}
      e^{-it_n^{j_0}\pt_x^4}f_n(x-x_n^{j_0})
      =\sum_{j=1}^{J} e^{i(t_n^j-t_n^{j_0})\pt_x^4}\phi^j(x-(x_n^j+x_n^{j_0})) 
      +e^{-it_n^{j_0}\pt_x^4}R_n^{J}(x-x_n^{j_0}). 
\end{align*}
Since $|x_n^{j_0}+x_n^j|\rightarrow\infty$ or $|t_n^{j_0}-t_n^j|\rightarrow\infty$ holds, we see from Lemma \ref{lem:t_n_x_n} that
\begin{align*}
      e^{i(t_n^j-t_n^{j_0})\pt_x^4}\phi^j(x-(x_n^j+x_n^{j_0})) 
      \rightharpoonup 0
\end{align*}
for all $j\in [1,J^{*}]\setminus\{j_0\}$. Hence, we obtain
\begin{align*}
      e^{-it_n^{j_0}\pt_x^4}R_n^{J}(x-x_n^{j_0})
      \rightharpoonup \phi^{j_0}(-x) \ \ \ \mbox{weakly in}\ H^2(\R).
\end{align*}
This implies that $\phi^{j_0}(-\cdot)\in\mathcal{V}(\{R_n^J\})$ for any finite $1\leq J\leq J^{*}$. Thus,  we have
\begin{align*}
  0<    \|\phi^{j_0}\|_{H^2}=\|\phi^{j_0}(-\cdot)\|_{H^2}
      \leq \eta(\{R_n^J\})
\end{align*} 
for any finite $1\leq J\leq J^{*}$. 
This contradicts with \eqref{eq:small_eta}.
On the other hand, if $x_n^j\equiv0$ then, arguing as in  \eqref{eq:pdevenpf}, one see that $\phi^j$ is an even function.
The property (v) is proved.
\qed \\

\begin{cor}\label{cor:energy_decomp.}
In the situation of Proposition \ref{prop:profile_decomp.}, for each finite $1\leq J\leq J^{*}$ we have
\begin{align*}
      E(f_n)=\sum_{j=1}^J E(e^{it_n^j\pt_x^4}\phi^j)+E(R_n^j)+o(1)
      \ \ \ \mbox{as}\ n\rightarrow\infty.
\end{align*}
\end{cor}

\noindent
{\bf Proof.} Corollary \ref{cor:energy_decomp.} can be proved by using \eqref{eq:norm_decomp.}, \eqref{eq:orthogonality} and the dispersive estimate \eqref{eq:dispersive}. 
\qed

\section{Concentration compactness}\label{sec:concentration_compactness}

For $0\leq L\leq M(Q)^{\frac{2-s_c}{s_c}}E(Q)$, we define
\begin{equation*}
      C(L):=\sup\left\{
      \|u\|_{X(\R)}\ \left|\ 
      \begin{split}
            &M(u)^{\frac{2-s_c}{s_c}}E(u)\leq L, \\
            &\|u_0\|_{L^2}^{\frac{2-s_c}{s_c}}\|u_0\|_{\dot{H}^2}
            <\|Q\|_{L^2}^{\frac{2-s_c}{s_c}}\|Q\|_{\dot{H}^2}
      \end{split}
      \right.
      \right\},
\end{equation*}
where the supremum is taken over all even solutions $u$ to \eqref{eq:1DNL4S} with initial data $u_0\in H^2(\R)$. From Proposition \ref{prop:stability}, $L\mapsto C(L)$ is non-decreasing and continuous at the every point in $[0,M(Q)^{\frac{2-s_c}{s_c}}E(Q))$ for which $C(L)$ is finite. From Proposition \ref{prop:small_data}, we see that $C(L)\lesssim L$ for all sufficiently small $L$. We set
\begin{align*}
      L^{*}:=\sup\{ L\in[0, M(Q)^{\frac{2-s_c}{s_c}}E(Q)]\ |\ C(L)<\infty \}.
\end{align*}
From the continuity of $C(L)$, we see that $C(L^{*})=\infty$. In the proof of Theorem \ref{thm:main}, the goal is to show $L^{*}=M(Q)^{\frac{2-s_c}{s_c}}E(Q)$. We will show this by a contradiction argument. In this section, we construct a special solution $u$ of \eqref{eq:1DNL4S} under the condition that $L^{*}<M(Q)^{\frac{2-s_c}{s_c}}E(Q)$. 

\begin{lem}[Concentration compactness]\label{lem:concentration_compactness}
Let $p>9$, and suppose $L^{*}<M(Q)^{\frac{2-s_c}{s_c}}E(Q)$. Let $\{u_n\}$ be a sequence of even solutions to \eqref{eq:1DNL4S} with initial data $u_{n,0}$ such that 
\begin{align*}
      &M(u_n)^{\frac{2-s_c}{s_c}}E(u_n)\leq L^{*}, \\
      &\|u_{n,0}\|_{L^2}^{\frac{2-s_c}{s_c}}\|u_{n,0}\|_{\dot{H}^2}
      <\|Q\|_{L^2}^{\frac{2-s_c}{s_c}}\|Q\|_{\dot{H}^2}
\end{align*}
for all $n$. Suppose that $M(u_n)^{\frac{2-s_c}{s_c}}E(u_n)\rightarrow L^{*}$ and there exists $\{t_n\}\subset\R$ such that
\begin{align*}
      \lim_{n\rightarrow\infty}\|u_n\|_{X((-\infty,t_n])}
      =\lim_{n\rightarrow\infty}\|u_n\|_{X([t_n,+\infty))}
      =\infty.
\end{align*}
Then there exists an even function $\phi\in H^2(\R)$ such that $u_n(t_n)$ has a subsequence that converges to $\phi$ in $H^2(\R)$. 
\end{lem}

\begin{rem}\label{rem:spatialcontrol}
The improvement due to the evenness assumption can be found in this lemma.
If we do not impose the assumption, we need to take spatial shift into account to obtain the convergence of
$\{ u_n(t_n)\}$ along a subsequence.
\end{rem}

\noindent
{\bf Proof.} Since $L^{*}<M(Q)^{\frac{2-s_c}{s_c}}E(Q)$, from Lemma \ref{lem:corecivity_1} we have
\begin{align*}
      \|u_n(t)\|_{L^2}^{\frac{2-s_c}{s_c}}\|u_n(t)\|_{\dot{H}^2}
      <\|Q\|_{L^2}^{\frac{2-s_c}{s_c}}\|Q\|_{\dot{H}^2}
\end{align*}
for all $t\in\R$. Thus, by using time translation symmetry, we may assume $t_n\equiv0$ without loss of generality. Then
\begin{align}\label{eq:blowup_X}
      \lim_{n\rightarrow\infty}\|u_n\|_{X((-\infty,0])}
      =\lim_{n\rightarrow\infty}\|u_n\|_{X([0,+\infty))}
      =\infty.
\end{align}
Applying Proposition \ref{prop:profile_decomp.} to $u_n(0)$, passing to a subsequence if necessary, for each $1\leq J\leq J^{*}$ we decompose
\begin{align*}
      u_n(0)=\sum_{j=1}^J e^{it_n^j\pt_x^4}\phi^j(x-x_n^j)+R_n^J,
\end{align*}
with
\begin{align*}
      &M(u_n)
      =\sum_{j=1}^J M(\phi^j)+M(R_n^J)+o(1),  \\ 
      &E(u_n)
      =\sum_{j=1}^J E(e^{it_n^j\pt_x^4}\phi^j)+E(R_n^J)+o(1), \\ 
      &\|u_n(0)\|_{\dot{H}^2}^2
      =\sum_{j=1}^J \|\phi^j\|_{\dot{H}^2}^2+\|R_n^J\|_{\dot{H}^2}^2+o(1), 
\end{align*}
where $o(1)\rightarrow0$ as $n\rightarrow\infty$, for each finite $1\leq J\leq J^{*}$. From these,
\begin{align}
      &\sum_{j=1}^J M(\phi^j)^{\frac{2-s_c}{s_c}}E(e^{it_n^j\pt_x^4}\phi^j)
      +M(R_n^J)^{\frac{2-s_c}{s_c}}E(R_n^J)
      \leq M(u_n)^{\frac{2-s_c}{s_c}}E(u_n)+o(1), \label{eq:sum_j_L} \\
      &\sum_{j=1}^J \|\phi^j\|_{L^2}^{\frac{2-s_c}{s_c}}\|\phi^j\|_{\dot{H}^2}
      +\|R_n^J\|_{L^2}^{\frac{2-s_c}{s_c}}\|R_n^J\|_{\dot{H}^2}
      \leq \|u_n(0)\|_{L^2}^{\frac{2-s_c}{s_c}}\|u_n(0)\|_{\dot{H}^2}+o(1) \label{eq:sum_j_scale_product}.
\end{align}
We claim that for each finite $1\leq J\leq J^{*}$,
\begin{align}
      &\inf_{1\leq j\leq J}\limsup_{n\rightarrow\infty}E(e^{it_n^j\pt_x^4}\phi^j)>0,\ \ \ 
      \limsup_{n\rightarrow\infty}E(R_n^J)
      \geq 0, \label{eq:inf_j_E} \\
      &\sup_{1\leq j\leq J}\limsup_{n\rightarrow\infty}
      M(\phi^j)^{\frac{2-s_c}{s_c}}E(e^{it_n^j\pt_x^4}\phi^j),\
      \limsup_{n\rightarrow\infty}M(R_n^J)^{\frac{2-s_c}{s_c}}E(R_n^J)
      \leq L^{*}, \label{eq:sup_j_L} \\
      &\sup_{1\leq j\leq J}\|\phi^j\|_{L^2}^{\frac{2-s_c}{s_c}}\|\phi^j\|_{\dot{H}^2},\
      \limsup_{n\rightarrow\infty}\|R_n^J\|_{L^2}^{\frac{2-s_c}{s_c}}\|R_n^J\|_{\dot{H}^2}
      < \|Q\|_{L^2}^{\frac{2-s_c}{s_c}}\|Q\|_{\dot{H}^2}. \label{eq:sup_j_scale_product}
\end{align}
Indeed, from \eqref{eq:sum_j_scale_product} we have
\begin{align}\label{eq:scale_product_profile}
      \sup_{1\leq j\leq J}\|\phi^j\|_{L^2}^{\frac{2-s_c}{s_c}}\|\phi^j\|_{\dot{H}^2},\
      \limsup_{n\rightarrow\infty}\|R_n^J\|_{L^2}^{\frac{2-s_c}{s_c}}\|R_n^J\|_{\dot{H}^2}
      \leq&\ \limsup_{n\rightarrow\infty}
      \|u_n(0)\|_{L^2}^{\frac{2-s_c}{s_c}}\|u_n(0)\|_{\dot{H}^2} \notag \\
      \leq&\ \|Q\|_{L^2}^{\frac{2-s_c}{s_c}}\|Q\|_{\dot{H}^2}.
\end{align}
Hence, by Lemma \ref{lem:corecivity_2} and $\phi^j\neq0$,  
\begin{align*}
      \inf_{1\leq j\leq J}\limsup_{n\rightarrow\infty}E(e^{it_n^j\pt_x^4}\phi^j)>0,\ \ \
      \limsup_{n\rightarrow\infty}E(R_n^J)
      \geq0.
\end{align*}
Thus, from \eqref{eq:sum_j_L} we have
\begin{align}
      &\sup_{1\leq j\leq J}\limsup_{n\rightarrow\infty}
      M(\phi^j)^{\frac{2-s_c}{s_c}}E(e^{it_n^j\pt_x^4}\phi^j),\
      \limsup_{n\rightarrow\infty}M(R_n^J)^{\frac{2-s_c}{s_c}}E(R_n^J) \notag \\
      \leq&\ \limsup_{n\rightarrow\infty}M(u_n)^{\frac{2-s_c}{s_c}}E(u_n)
      \leq L^{*}. \label{eq:L_linear_profile}
\end{align}
Since $L^{*}<M(Q)^{\frac{2-s_c}{s_c}}E(Q)$, By using \eqref{eq:scale_product_profile}, \eqref{eq:L_linear_profile} and the sharp Gagliardo--Nirenberg inequality, we have
\begin{align*}
      \sup_{1\leq j\leq J}\|\phi^j\|_{L^2}^{\frac{2-s_c}{s_c}}\|\phi^j\|_{\dot{H}^2},\
      \limsup_{n\rightarrow\infty}\|R_n^J\|_{L^2}^{\frac{2-s_c}{s_c}}\|R_n^J\|_{\dot{H}^2}
      <\|Q\|_{L^2}^{\frac{2-s_c}{s_c}}\|Q\|_{\dot{H}^2}.
\end{align*}
\indent
Now we consider the following three cases: $J^{*}=0$, $J^{*}=1$ and $J^{*}\geq2$. We will show that the first case and the third case lead to contradiction. In the second case, we will get desired $\phi\in H^2(\R)$. \\

\indent
{\bf Case 1:} $J^{*}=0$. In this case, we have
\begin{align*}
      \limsup_{n\rightarrow\infty}\|e^{-it\pt_x^4}u_n(0)\|_{X(\R)}
      =\limsup_{n\rightarrow\infty}\|e^{-it\pt_x^4}R_n^0\|_{X(\R)}
      =0.
\end{align*}
Hence, for $n$ sufficiently large, we obtain
\begin{align*}
      \|e^{-it\pt_x^4}u_n(0)\|_{X(\R)}\leq \delta_{sd}.
\end{align*}
Thus, by using Proposition \ref{prop:small_data}, we have
\begin{align*}
      \|u_n\|_{X(\R)}\leq2\|e^{-i\pt_x^4}u_n(0)\|_{X(\R)}\leq2\delta_{sd}.
\end{align*}
This contradicts \eqref{eq:blowup_X} and so $J^{*}=0$ does not occur. \\

\indent
{\bf Case 2:} $J^{*}=1$. In this case, there exists a single profile in the decomposition and so we can write
\begin{align*}
      u_n(0)=e^{it_n\pt_x^4}\phi^1(x-x_n)+R_n^1.
\end{align*}
By the improvement of the profile decomposition by the evenness assumption ((v) of Propoition \ref{prop:profile_decomp.}), one sees that 
$|x_n|$ is bounded and hence $x_n \equiv 0$.
Further, $\phi^1$ is an even function.

Let us prove $t_n\equiv0$. 
If $t_n\rightarrow-\infty$, then by the Strichartz estimate and the monotone convergence theorem,
\begin{align*}
      \|e^{-i(t-t_n)\pt_x^4}\phi^1\|_{X([0,+\infty))}
      = \|e^{-it\pt_x^4}\phi^1\|_{X([-t_n,+\infty))}
      \rightarrow 0
\end{align*}
as $n\rightarrow\infty$. From this and \eqref{eq:smallness}, we have $\limsup_{n\rightarrow\infty}\|e^{-it\pt_x^4}u_n(0)\|_{X([0,+\infty))}=0$. Hence, by Proposition \ref{prop:small_data}, for $n$ sufficiently large, we obtain 
\begin{align*}
      \|u_n\|_{X([0,+\infty))}\leq2\delta_{sd}.
\end{align*}
This contradicts \eqref{eq:blowup_X} and so $t_n\rightarrow-\infty$ does not occur. By the similar argument, $t_n\rightarrow+\infty$ does not occur. 
Thus, we have $t_n\equiv0$ and hence
\begin{align*}
      u_n(0)=\phi^1+R_n^1.
\end{align*}
Then we claim $M(\phi^1)^{\frac{2-s_c}{s_c}}E(\phi^1)=L^{*}$. To prove this, suppose $M(\phi^1)^{\frac{2-s_c}{s_c}}E(\phi^1)<L^{*}$ and let $v$ be a solution to \eqref{eq:1DNL4S} with $v(0)=\phi^1$. Then by the definition of $L^{*}$, we have $\|v\|_{X(\R)}<\infty$. Let
\begin{align*}
      \tilde{u}_n(t):=v(t)+e^{-it\pt_x^4}R_n^1
\end{align*}
and 
\begin{align*}
      e_n:=&\ i\pt_t \tilde{u}_n-\pt_x^4 \tilde{u}_n+|\tilde{u}_n|^{p-1}\tilde{u}_n \notag \\
      =&\ |\tilde{u}_n|^{p-1}\tilde{u}_n-|v(t)|^{p-1}v(t).
\end{align*}
Then $\tilde{u}_n(0)=u_n(0)$ and by Lemma \ref{lem:nonlinear_est}, we have
\begin{align*}
      \|e_n\|_{N(\R)}
      \lesssim&\ \left(\|v\|_{X(\R)}+\|e^{-it\pt_x^4}R_n^1\|_{X(\R)}\right)^{p-1}
      \|e^{-it\pt_x^4}R_n^1\|_{X(\R)}.
\end{align*}
Hence, using $\|v\|_{X(\R)}<\infty$ and $\limsup_{n\rightarrow\infty}\|e^{-it\pt_x^4}R_n\|_{X(\R)}=0$, we obtain 
\begin{align*}
      \limsup_{n\rightarrow\infty}\|e_n\|_{N(\R)}=0.
\end{align*}
Thus, by using Proposition \ref{prop:stability}, 
\begin{align}\label{eq:X_tilde_u_n}
      \|u_n\|_{X(\R)}
      \leq \|\tilde{u}_n\|_{X(\R)}+o(1)
\end{align}
as $n\rightarrow\infty$. Since $\|\tilde{u}_n\|_{X(\R)}\leq \|v\|_{X(\R)}+\|e^{-it\pt_x^4}R_n^1\|_{X(\R)}<\infty$ for sufficiently large $n$, \eqref{eq:X_tilde_u_n} contradicts \eqref{eq:blowup_X} and so $M(\phi^1)^{\frac{2-s_c}{s_c}}E(\phi^1)=L^{*}$. In particular, we have $\limsup_{n\rightarrow\infty}\|R_n\|_{H^2}=0$. Therefore, 
\begin{align*}
      u_n(0)=\phi^1+o(1)
\end{align*}
in $H^2(\R)$. 
This completes the proof of Lemma \ref{lem:concentration_compactness} if we have proved that $J^{*}\geq2$ does not occur.
\\

\indent
{\bf Case 3:} $J^{*}\geq2$. In this case, there exist more than one profiles, and from \eqref{eq:inf_j_E} and \eqref{eq:sup_j_L} there exists $\epsilon>0$ such that 
\begin{align}\label{eq:L_of_v^j}
      0<\limsup_{n\rightarrow\infty}M(\phi^j)^{\frac{2-s_c}{s_c}}E(e^{it_n^j\pt_x^4}\phi^j)
      \leq L^{*}-\epsilon
\end{align}
for all finite $1\leq j\leq J^{*}$. For each finite $1\leq j\leq J^{*}$, passing to a subsequence, we may assume that $t_n^j$ has a limit $t_0^j\in[-\infty,+\infty]$ as $n\rightarrow\infty$. By using Lemma \ref{prop:small_data}, we can find an unique solution $v^j$ of \eqref{eq:1DNL4S} defined on a neighborhood of $-t_0^j$ such that
\begin{align}\label{eq:H^2_v^j}
      \|v^j(-t_n^j)-e^{it_n^j\pt_x^4}\phi^j\|_{H^2}\rightarrow0\ \ \ \mbox{as}\ n\rightarrow\infty.
\end{align}
From \eqref{eq:sup_j_L}, \eqref{eq:sup_j_scale_product} and \eqref{eq:H^2_v^j} we have
\begin{align*}
      M(v^j)^{\frac{2-s_c}{s_c}}E(v^j)\leq L^{*}<M(Q)^{\frac{2-s_c}{s_c}}E(Q)
\end{align*}
and
\begin{align*}
      \|v^j(t)\|_{L^2}^{\frac{2-s_c}{s_c}}\|v^j(t)\|_{\dot{H}^2}
      <\|Q\|_{L^2}^{\frac{2-s_c}{s_c}}\|Q\|_{\dot{H}^2}
\end{align*}
for all $t\in\R$. Hence, from Lemma \ref{lem:corecivity_1}, we see that $v^j$ exists globally in time. Let $v_n^j:=v^j(t-t_n^j,x-x_n^j)$, and for each finite $1\leq J\leq J^{*}$ we define
\begin{align*}
      \tilde{u}_n^J:=\sum_{j=1}^J v_n^j+e^{it\pt_x^4}R_n^J.
\end{align*}
Then 
\begin{align*}
      \|\tilde{u}_n^J\|_{X(\R)}
      \leq&\ \sum_{j=1}^J \|v^j\|_{X(\R)}
      +C_J\sum_{1\leq j\neq k\leq J}\||v_n^j|
      |v_n^k|^{r-1}\|_{L_t^{q_1/r}L_x^1(\R\times\R)}
      +\|e^{it\pt_x^4}R_n^J\|_{X(\R)}, 
\end{align*}
where $C_J>0$ is a constant depend only on $J$. Since $C(L)$ is sublinear around $0$ and bounded on $[0,L^{*}-\epsilon]$, from \eqref{eq:L_of_v^j} we obtain
\begin{align}\label{eq:X_sum_v_n^j}
      \sum_{j=1}^J \|v^j\|_{X(\R)}
      \leq&\ \sum_{j=1}^J C\left(M(v^j)^{\frac{2-s_c}{s_c}}E(v^j)\right) \notag \\
      \lesssim&_{L^{*},\epsilon}\ \sum_{j=1}^J M(v^j)^{\frac{2-s_c}{s_c}}E(v^j) \notag \\
      \lesssim&_{L^{*},\epsilon}\ 1.
\end{align}
For each finite $1\leq J\leq J^{*}$, using \eqref{eq:orthogonality}, we can take $n$ sufficiently large so that
\begin{align*}
      \sum_{1\leq j\neq k\leq J}\||v_n^j||v_n^k|^{r-1}\|_{L_t^{q_1/r}L_x^1(\R\times\R)}
      \leq \frac{1}{C_J J^2}.
\end{align*}
Therefore, by using \eqref{eq:smallness}, we obtain
\begin{align}\label{eq:X_tilde_u_n^J}
      \lim_{J\rightarrow J^{*}}\limsup_{n\rightarrow\infty}\|\tilde{u}_n^J\|_{X(\R)}
      \lesssim_{L^{*},\epsilon}1.
\end{align}
We define
\begin{align*}
      e_n^J:= i\pt_t \tilde{u}_n^J-\pt_x^4 \tilde{u}_n^J+|\tilde{u}_n^J|^{p-1}\tilde{u}_n^J 
      = F\left(\sum_{j=1}^J v_n^j+e^{it\pt_x^4}R_n^J\right)-\sum_{j=1}^J F(v_n^j), 
\end{align*}
where $F(z)=|z|^{p-1}z$ for $z\in\C$. Since
\begin{align*}
      &\left\|F\left(\sum_{j=1}^J v_n^j +e^{it\pt_x^4}R_n^J\right)
      -F\left(\sum_{j=1}^J v_n^j\right)\right\|_{N(\R)} \notag \\
      &\lesssim\ \left\| \left( \left|\sum_{j=1}^J v_n^j\right|^{p-1}+|e^{it\pt_x^4}R_n^J|^{p-1} \right)
      |e^{it\pt_x^4}R_n^J| \right\|_{N(\R)} \notag \\
      &\lesssim\ \left( \sum_{j=1}^J \|v^j\|_{X(\R)}^{p-1}+\|e^{it\pt_x^4}R_n^J\|_{X(\R)}^{p-1} \right)
      \|e^{it\pt_x^4}R_n^J\|_{X(\R)},
\end{align*}
from \eqref{eq:smallness} and \eqref{eq:X_sum_v_n^j}, we have
\begin{align*}
      \lim_{J\rightarrow J^{*}}\limsup_{n\rightarrow\infty}
      \left\|F\left(\sum_{j=1}^J v_n^j +e^{it\pt_x^4}R_n^J\right)
      -F\left(\sum_{j=1}^J v_n^j\right)\right\|_{N(\R)}
      =0.
\end{align*}
By using \eqref{eq:orthogonality}, we obtain
\begin{align*}
      \left\| F\left(\sum_{j=1}^J v_n^j\right)-\sum_{j=1}^J F(v_n^j) \right\|_{N(\R)}
      \lesssim C_J \sum_{1\leq j\neq k\leq J} \||v_n^j||v_n^k|^{p-1}\|_{N(\R)}
      =o_J(1).
\end{align*}
Therefore, we get
\begin{align}\label{eq:N_e_n^J}
      \lim_{J\rightarrow J^{*}}\limsup_{n\rightarrow\infty}\|e_n^J\|_{N(\R)}
      =0.
\end{align}
From \eqref{eq:H^2_v^j} we have
\begin{align}\label{eq:tilde_u_n^J(0)}
      \tilde{u}_n^J(0)
      =\sum_{j=1}^J v^j(-t_n^j,x-x_n^j) +R_n^J
      =u_n(0)+o_J(1)
\end{align}
as $n\rightarrow\infty$ in $H^2(\R)$. From \eqref{eq:X_tilde_u_n^J}, \eqref{eq:N_e_n^J} and \eqref{eq:tilde_u_n^J(0)}, using Proposition \ref{prop:stability}, we have
\begin{align*}
      \limsup_{n\rightarrow\infty}\|u_n\|_{X(\R)}
      \leq \lim_{J\rightarrow J^{*}}\limsup_{n\rightarrow\infty}\|\tilde{u}_n^J\|_{X(\R)}
      \lesssim_{L^{*},\epsilon} 1.
\end{align*}
This contradicts \eqref{eq:blowup_X} and so $J^{*}\geq2$ does not occur.
\qed \\

\begin{prop}[Existence and precompactness of critical solution]\label{prop:critical_sol}
Let $p>9$. Suppose $L^{*}<M(Q)^{\frac{2-s_c}{s_c}}E(Q)$. Then there exists $\phi\in H^2(\R)$ such that
\begin{align}\label{eq:L_of_phi}
      &M(\phi)^{\frac{2-s_c}{s_c}}E(\phi)=L^{*}, \\
      &\|\phi\|_{L^2}^{\frac{2-s_c}{s_c}}\|\phi\|_{\dot{H}^2}
      <\|Q\|_{L^2}^{\frac{2-s_c}{s_c}}\|Q\|_{\dot{H}^2}. \notag
\end{align}
Let $u$ be a solution to \eqref{eq:1DNL4S} with initial data $\phi$. Then
\begin{align}\label{eq:X_of_u}
      \|u\|_{X((-\infty,0])}=\|u\|_{X([0,+\infty))}=\infty.
\end{align}
Moreover, $\{ u(t)\ |\ t\in\R \}$ is precompact in $H^2(\R)$.
\end{prop}

\begin{rem}\label{rem:spatialcontrol2}
Because of the evenness assumption, the precompactness
of the orbit holds without spatial translation.
As is clear from the proof, this is due to the improvement in Lemma  \ref{lem:concentration_compactness} discussed in 
Remark \ref{rem:spatialcontrol}.
The property is nothing but the long-time localization property of the critical element.
\end{rem}

\noindent
{\bf Proof.} Since $C(L):[0,L^{*})\rightarrow[0,\infty)$ is continuous, non-decreasing and $C(L^{*})=\infty$, we can take $L_n$ such that $C(L_n)=n$ and $L_n\nearrow L^{*}$ as $n\rightarrow\infty$. Then by the definition of $C(L_n)$, there exists a sequence of even solutions $\{u_n\}$ to \eqref{eq:1DNL4S} with a initial data $u_{n,0}$ such that 
\begin{align}\label{eq:L_n}
      L_{n-1}\leq M(u_n)^{\frac{2-s_c}{s_c}}E(u_n)< L_n, \ \ \ 
      \|u_{n,0}\|_{L^2}^{\frac{2-s_c}{s_c}}\|u_{n,0}\|_{\dot{H}^2}
      <\|Q\|_{L^2}^{\frac{2-s_c}{s_c}}\|Q\|_{\dot{H}^2}
\end{align}
and $\|u_n\|_{X(\R)}\geq n-\frac{1}{2}$. Let $f(t):=\|u_n\|_{X((-\infty,t])}$. Then $f(-\infty)=0$ and $f(+\infty)\geq n-\frac{1}{2}$. Since $f(t)$ is continuous, we find $t_n\in\R$ such that $\frac{1}{4}(n-\frac{1}{2})\leq f(t_n)\leq \frac{1}{2}(n-\frac{1}{2})$. Then
\begin{align*}
      &\|u_n\|_{X((-\infty,t_n])}=f(t_n)\geq \frac{1}{4}(n-\frac{1}{2})\rightarrow\infty, \notag \\
      &\|u_n\|_{X([t_n,+\infty)}=\|u_n\|_{X(\R)}-f(t_n)\geq \frac{1}{2}(n-\frac{1}{2})\rightarrow\infty.
\end{align*}
Moreover, by Lemma \ref{lem:corecivity_1}, 
\begin{align}\label{eq:u_n(t_n)}
      \|u_n(t_n)\|_{L^2}^{\frac{2-s_c}{s_c}}\|u_n(t_n)\|_{\dot{H}^2}
      <\|Q\|_{L^2}^{\frac{2-s_c}{s_c}}\|Q\|_{\dot{H}^2}.
\end{align}
Applying Lemma \ref{lem:concentration_compactness}, there exists an even function $\phi\in H^2(\R)$ such that $u_n(t_n)\rightarrow\phi$ in $H^2(\R)$ as $n\rightarrow\infty$. From this strong convergence, \eqref{eq:L_n} and \eqref{eq:u_n(t_n)}, we have \eqref{eq:L_of_phi}. Let $u$ be a solution to \eqref{eq:1DNL4S} with initial data $\phi$. Then by using Propsition \ref{prop:stability}, we obtain \eqref{eq:X_of_u}. \\
\indent
Finally, we prove that $\{ u(t)\ |\ t\in\R \}$ is precompact in $H^2(\R)$. Take a sequence $\{t_n^{\prime}\}\subset\R$. Then from \eqref{eq:X_of_u} we have
\begin{align*}
      \|u\|_{X((-\infty,t_n^{\prime}])}=\|u\|_{X([t_n^{\prime},+\infty))}
      =\infty.
\end{align*}
Applying Lemma \ref{lem:concentration_compactness} again, there exists $\phi^{\prime}\in H^2(\R)$ such that $u(t_n^{\prime})$ has a subsequence that converges $\phi^{\prime}$ in $H^2(\R)$. This completes the proof of Proposition \ref{prop:critical_sol}.
\qed

\section{Rigidity}\label{sec:rigidity}

In this section, we will complete the proof of Theorem \ref{thm:main}. We first prepare the following Lemma, in which precompactness in $H^2(\R)$ implies tightness in $H^2(\R)$.

\begin{lem}[Tightness in $H^2$]\label{lem:tightness}
Let $u$ be a global solution to \eqref{eq:1DNL4S} such that $\{u(t)\ |\ t\in\R\}$ is precompact in $H^2(\R)$. Then for all $\eta>0$, there exists $R_0>0$ such that for all $R\geq R_0$ and for all $t\in\R$,
\begin{align}\label{eq:tightness}
      \int_{|x|\geq R}|u(t,x)|^2+|\pt_x^2 u(t,x)|^2 dx
      <\eta.
\end{align}
\end{lem}

\noindent
{\bf Proof.} Suppose that Lemma \ref{lem:tightness} does not hold. Then there exist $\eta_0>0$, $R_n>0$ and $t_n\in\R$ such that $R_n\rightarrow\infty$ and
\begin{align}\label{eq:>eta_0}
      \int_{|x|\geq R_n}|u(t_n,x)|^2+|\pt_x^2 u(t_n,x)|^2 dx \geq \eta_0
\end{align}
for all $n$. Since $\{ u(t)\ |\ t\in\R \}$ is precompact in $H^2(\R)$, there exists $\phi\in H^2(\R)$ such that $u(t_n)$ has a subsequence that converges $\phi$ in $H^2(\R)$. In particular, there exists $n_1$ such that
\begin{align}\label{eq:H^2_u-phi}
      \|u(t_n)-\phi\|_{H^2}<\frac{\eta_0}{4}
\end{align}
for all $n\geq n_1$. On the other hand, since $R_n\rightarrow\infty$, there exists $n_2$ such that
\begin{align}\label{eq:H^2(|x|>R_n)_phi}
      \int_{|x|\geq R_n} |\phi(x)|^2+|\pt_x^2 \phi(x)|^2 dx<\frac{\eta_0}{4}
\end{align}
for all $n\geq n_2$. By the triangle inequality, \eqref{eq:H^2_u-phi} and \eqref{eq:H^2(|x|>R_n)_phi} contradict \eqref{eq:>eta_0} when $n\geq\max\{n_1,n_2\}$.
\qed \\

Next, we show that there is no nontrivial solution which has a precompcat orbit in $H^2(\R)$ and satisfies \eqref{eq:below_ground_state}. To prove this, we use the localized virial identity associated with fourth-order nonlinear Schr\"{o}dinger equation \eqref{eq:1DNL4S}.

\begin{prop}[Rigidity]\label{prop:rigidity}
Let $p>9$. Let $u$ be a solution to \eqref{eq:1DNL4S} with initial data $u_0$ such that 
\begin{align*}
      M(u)^{\frac{2-s_c}{s_c}}E(u)<&M(Q)^{\frac{2-s_c}{s_c}}E(Q), \notag \\
      \|u_0\|_{L^2}^{\frac{2-s_c}{s_c}}\|u_0\|_{\dot{H}^2}
      <&\|Q\|_{L^2}^{\frac{2-s_c}{s_c}}\|Q\|_{\dot{H}^2}
\end{align*}
and $\{ u(t)\ |\ t\in\R \}$ is precompact in $H^2(\R)$. Then $u_0=0$.
\end{prop}

\noindent
{\bf Proof.} Suppose $u_0\neq0$. From Lemma \ref{lem:corecivity_1}, $u$ exists globally in time. In particular, there exists $\delta^{\prime}\in(0,1)$ such that 
\begin{align}\label{eq:u_delta'_Q}
      \|u(t)\|_{L^2}^{\frac{2-s_c}{s_c}}\|u(t)\|_{\dot{H}^2}
      <(1-\delta^{\prime})\|Q\|_{L^2}^{\frac{2-s_c}{s_c}}\|Q\|_{\dot{H}^2}
\end{align}
for all $t\in\R$. This implies that there exists a constant $c_1=c_1(\delta^{\prime},Q)>0$ such that $\|u(t)\|_{H^2}\leq c_1$ for all $t\in\R$. Let $\psi\in C_0^{\infty}([0,\infty))$ be a non-increasing function and such that $\psi(s)=1$ for $0\leq s\leq1$ and $\psi(r)=0$ for $s\geq2$. For each $R>0$, we define 
\begin{align*}
      \Psi_R(x):=R^2\int_{0}^{\frac{|x|}{R}}\int_0^r \psi(s) dsdr
\end{align*} 
and
\begin{align*}
      M_R(t):=2\im\int_{\R}\overline{u}(t,x)\pt_x\Psi_R(x)\pt_x u(t,x) dx.
\end{align*}
Then by the H\"{o}lder inequality, 
\begin{align*}
|M_R(t)|\lesssim R\|u(t)\|_{L^2}\|\pt_x u(t)\|_{L^2}\leq R\|u(t)\|_{H^2}^2\leq c_1^2 R
\end{align*}
for all $t\in\R$. By using \eqref{eq:1DNL4S} and applying the integration by parts, we obtain
\begin{align*}
      \frac{d}{dt}M_R(t)
      =&\ 8\int_{\R}\pt_x^2\Psi_R|\pt_x^2 u|^2 dx-6\int_{\R}\pt_x^4\Psi_R|\pt_x u|^2 dx 
      +\int_{\R}\pt_x^6\Psi_R|u|^2 dx \notag \\
      &\ -\frac{2(p-1)}{p+1}\int_{\R}\pt_x^2\Psi_R|u|^{p+1} dx \notag \\
      =&\ 4K(u(t))+A_R(t),
\end{align*}
where
\begin{align*}
      A_R(t)=&\ 8\int_{\R}(\pt_x^2\Psi_R-1)|\pt_x^2 u|^2 dx
      -6\int_{\R}\pt_x^4\Psi_R|\pt_x u|^2 dx 
      +\int_{\R}\pt_x^6\Psi_R|u|^2 dx \notag \\
      &\ -\frac{2(p-1)}{p+1}\int_{\R}(\pt_x^2\Psi_R-1)|u|^{p+1} dx.
\end{align*}
Since $\pt_x^2\Psi_R(x)=\psi(\frac{|x|}{R})$, by using Lemma \ref{lem:tightness}, for any $\eta>0$ we can take $R>0$ sufficiently large so that 
\begin{align}\label{eq:A_R_1}
      \left|\ 8\int_{\R}(\pt_x^2\Psi_R-1)|\pt_x^2 u(t,x)|^2 dx\ \right|
      \leq8\int_{|x|\geq R}|\pt_x^2 u(t,x)|^2 dx
      <\eta
\end{align}
for all $t\in\R$. Since $\pt_x^4\Psi_R(x)=R^{-2}\psi^{\prime\prime}(\frac{|x|}{R})$ and $\pt_x^6\Psi_R(x)=R^{-4}\psi^{(4)}(\frac{|x|}{R})$, we have
\begin{align}\label{eq:A_R_2}
      \left|\ -6\int_{\R}\pt_x^4\Psi_R|\pt_x u(t,x)|^2 dx
      +\int_{\R}\pt_x^6\Psi_R|u(t,x)|^2 dx\ \right|
      \lesssim&\ (R^{-2}+R^{-4})\|u(t)\|_{H^2}^2 \notag \\
      \leq&\ (R^{-2}+R^{-4})c_1^2,
\end{align}
for all $t\in\R$. Since $\pt_x^2\Psi_R(x)=\psi(\frac{|x|}{R})$ and so $\pt_x^2\Psi_R(x)-1\leq0$, we have
\begin{align}\label{eq:A_R_3}
      -\frac{2(p-1)}{p+1}\int_{\R}(\pt_x^2\Psi_R-1)|u(t,x)|^{p+1} dx\geq0.
\end{align}
From Lemma \ref{lem:corecivity_2}, there exists $\delta^{\prime\prime}\in(0,1)$ such that 
\begin{align*}
      K(u(t))\geq 2\delta^{\prime\prime}\|u(t)\|_{\dot{H}^2}^2
      \geq 4\delta^{\prime\prime}E(u)
\end{align*}
for all $t\in\R$. Moreover,
\begin{align*}
      E(u(t))\geq \frac{p-9}{2(p-1)}\|u(t)\|_{\dot{H}^2}^2
\end{align*}
for all $t\in\R$. Thus, if $u_0\neq0$, then $E(u)$ is positive. Therefore, from \eqref{eq:A_R_1}, \eqref{eq:A_R_2} and \eqref{eq:A_R_3}, we can take $R>0$ sufficiently large so that $A_R(t)\geq-8\delta^{\prime\prime}E(u)$ and so
\begin{align*}
      \frac{d}{dt}M_R(t)\geq 8\delta^{\prime\prime}E(u)
      \geq \frac{4(p-9)}{p-1}\|u_0\|_{\dot{H}^2}^2.
\end{align*}
Hence, for any time interval $I$,
\begin{align*}
      \delta^{\prime\prime}\|u_0\|_{\dot{H}^2}^2|I|
      \lesssim \sup_{t\in\R}|M_R(t)|
      \lesssim c_1^2 R,
\end{align*} 
which contradicts $u_0\neq0$ for sufficiently large $|I|$. 
\qed \\

\noindent
{\bf Proof of Theorem \ref{thm:main}.} To prove Theorem \ref{thm:main}, it suffices to show $L^{*}\geq M(Q)^{\frac{2-s_c}{s_c}}E(Q)$. Indeed, if $L^{*}\geq M(Q)^{\frac{2-s_c}{s_c}}E(Q)$, then $C(L)<\infty$ for all $L<M(Q)^{\frac{2-s_c}{s_c}}E(Q)$. Thus, by the definition of $C(L)$, all even solutions $u$ of \eqref{eq:1DNL4S} such that \eqref{eq:below_ground_state} satisfy $\|u\|_{X(\R)}<\infty$. Hence, by using Proposition \ref{prop:H^2-scattering}, u scatters forward and backward in time. \\
\indent
Now we will prove $L^{*}\geq M(Q)^{\frac{2-s_c}{s_c}}E(Q)$. Suppose that $L^{*}<M(Q)^{\frac{2-s_c}{s_c}}E(Q)$, then from Proposition \ref{prop:critical_sol}, there exists a global solution $u$ of \eqref{eq:1DNL4S} with an initial data $u_0$ such that 
\begin{align*}
      M(u)^{\frac{2-s_c}{s_c}}E(u)<&\ M(Q)^{\frac{2-s_c}{s_c}}E(Q), \notag \\
      \|u_0\|_{L^2}^{\frac{2-s_c}{s_c}}\|u_0\|_{\dot{H}^2}
      <&\ \|Q\|_{L^2}^{\frac{2-s_c}{s_c}}\|Q\|_{\dot{H}^2}
\end{align*}
and 
\begin{align}\label{eq:blowup_u}
      \|u\|_{X((-\infty,0])}=\|u\|_{X([0,+\infty))}=\infty.
\end{align}
Moreover, $\{ u(t)\ |\ t\in\R \}$ is precompact in $H^2(\R)$. By using Proposition \ref{prop:rigidity}, $u_0=0$. Hence, $u\equiv0$ which contradicts \eqref{eq:blowup_u}. Thus, $L^{*}\geq M(Q)^{\frac{2-s_c}{s_c}}E(Q)$.
\qed

\subsection*{Acknowledgement}
K. K. was supported by JSPS KAKENHI Grant Number 23K13003.
S. M. was supported by JSPS KAKENHI Grant Numbers  	21H00991 and 21H00993.


\end{document}